\documentclass[11pt]{article}
\pdfoutput=1
\usepackage{amssymb}
\usepackage{amsfonts}
\usepackage{amsmath}
\usepackage{amscd}
\usepackage[T1]{fontenc}
\usepackage{graphicx}

\newcommand{\R}{\mathbb{R}}
\newcommand{\C}{\mathbb{C}}

\newcommand{\Z}{\mathbb{Z}}

\newcommand{\Gl}{\mathop{\mathrm{Gl}}\nolimits}

\newcommand{\SSl}{\mathop{\mathrm{Sl}}\nolimits}

\newcommand{\dee}{\mathop{\! \, \mathrm{d} \!}\nolimits}
\newcommand{\comp}{\raisebox{0pt}{$\scriptstyle\circ \, $}}

\newcommand{\setrule}{\, \mathop{\rule[-4pt]{.5pt}{13pt}\, }\nolimits}
\newcommand{\bigspace}{\medskip\par\noindent}
\newcommand{\smallspace}{\smallskip\par\noindent}
\newcommand{\smallrowspace}{\rule{0pt}{12pt}}
\newcommand{\rowspace}{\rule{0pt}{16pt}}
\newcommand{\largerowspace}{\rule{0pt}{20pt}}

\newcommand{\onehalf}{\mbox{$\frac{\scriptstyle 1}{\scriptstyle 2}\,$}}

\newcommand{\dbydt}{\mbox{${\displaystyle \frac{\dee}{\dee t}}
\rule[-10pt]{.5pt}{25pt} \raisebox{-10pt}{$\, {\scriptstyle t=0}$}$}}
\newcommand{\ttfrac}[2]{\mbox{$\frac{{\scriptstyle #1}}{{\scriptstyle #2}}$}}
\newcommand{\poisbr}[3]{\mbox{${\{ #1 , #2 \}} _{{\scriptstyle #3}}$}}
\newcommand{\vvee}{\mbox{\tiny $\vee $}}
\allowdisplaybreaks

\begin{document}

\begin{center}
{\Large \bf Classical and Quantum Spherical Pendulum} 
\end{center}
\hspace{1in} {Richard Cushman\footnote{Department of Mathematics and Statistics, 
University of Calgary, \\ email: rcushman(at)ucalgary.ca} and\ J\k{e}drzej \'{S}niatycki}\footnote{Departments of Mathematics and Statistics, University of Calgary and University of Victoria, \\
email:sniatyck(at)ucalgary.ca} 

\begin{abstract}
This paper extends the Bohr-Sommerfeld quantization of the spherical
pendulum to a full quantum theory. This the first application of geometric quantization to a classical system
with monodromy.
\end{abstract}

\section{Introduction}

This paper is a part of a program to extend Bohr-Sommerfeld conditions to a
full quantum theory of completely integrable systems. In earlier
publications, we discussed several examples of systems with one and two
degrees of freedom \cite{cushman-sniatycki13}, \cite{cushman-sniatycki14}.
In the present paper, we concentrate our attention on the spherical
pendulum. Our aim is to understand quantum monodromy of this system 
\cite{cushman-duistermaat}. More precisely, we want to understand how the
classical monodromy emerges from the quantum theory of the spherical
pendulum.

The sperical pendulum has global actions $A_{1},A_{2}$, where $A_{1}$ is a
continuous function of the energy $H$ and the angular momentum $L$, and $%
A_{2}=L$. The Bohr-Sommerfeld conditions $A_{1}=2\pi n\hbar $ and $%
A_{2}=2\pi m\hbar ,$ where $n$ and $m$ are integers and $\hbar $ is Planck's
constant divided by $2\pi $, give quantization of the action functions $%
A_{1},A_{2}$, and determine the joint spectrum of the operators ${\mathbf{Q}}_{H}$ and $%
{\mathbf{Q}}_{L}$ corresponding to $H$ and $L$, respectively. Moreover, the pairs $%
(n,m) $ of integers, for which the Bohr-Sommerfeld conditions are satisfied,
label the basic states of a basis $\{\sigma _{n,m}\}$ of the space of quantum
states $\mathfrak{H}$ of the system. In geometric quantization, each $\sigma _{nm}$ is a distribution 
section (a generalized section) of the
prequantization line bundle with support satisfying Bohr-Sommerfeld
conditions with integers $(n,m)$ \cite{sniatycki80}. This implies that the
basis $\{\sigma _{n,m}\}$ of $\mathfrak{H}$ has the structure of a local
lattice. This observation lead Cushman and Duistermaat to the notion of
quantum mondromy \cite{cushman-duistermaat}. In the following, we refer to
the basis $\{\sigma _{n,m}\}$ as the Bohr-Sommerfeld basis.

The Bohr-Sommerfeld theory, as outlined above, does not lead to operators on
the space $\mathfrak{H}$ of quantum states that are not diagonal in the
basis $\{\sigma _{n,m}\}$. Though such operators exist, the theory does not
relate them to classical functions. In his 1925 paper \cite{heisenberg},
Heisenberg emphasized the importance of operators that provide transitions
between different quantum states. Since Bohr-Sommerfeld approach did not
provide transition operators, it was abandoned in favour of the matrix
mechanics of Born, Jordan and Heisnberg \cite{born-jordan}, \cite%
{born-jordan-heisenberg} and the wave mechanics of Schr\"{o}dinger \cite{schrodinger}. Dirac incorporated both approaches in his \textit{Principles of Quantum Mechanics }published in 1930 \cite{dirac 30}.

In \cite{cushman-sniatycki14}, we showed that, if the Bohr-Sommerfeld basis
of the space $\mathfrak{H}$ of quantum states of a completely integrable
system has a structure of a global lattice with boundary and the
action-angle variables are globally defined in the open dense set of regular
points of the energy-momentum map $\mathcal{EM}$, then the lowering
operators that take $\sigma _{n,m}$ to $\sigma _{n-1,m}$ can be interpreted
as quantization of ${\mathrm{e}}^{-i\varphi _{1}}$. Similarly, the operators that take $%
\sigma _{n,m}$ to $\sigma _{n,m-1\textrm{ }}$can be interpreted as
quantization of ${\mathrm{e}}^{-i\varphi _{2}}$. In this paper, we show that the
Bohr-Sommerfeld basis of the space of quantum states of the spherical
pendulum has the structure of the global lattice with boundary so that
shifting operators are well defined. However, due to the classical
monodromy, our action-angle variables fail to satisfy the defining equation 
$\omega = \dee A_{1}\wedge \dee \varphi _{1}+ \dee A_{2}\wedge \dee \varphi _{2} $
on $L^{-1}(0)$. Thus, Dirac's quantization conditions allow us to interperate the lowering 
operators as quantizations of ${\mathrm{e}}^{-i\varphi _{1}}$ and ${\mathrm{e}}^{-i\varphi _{2}}$ 
only in the complement of $L^{-1}(0)$ in $\mathcal{EM}^{-1}(\mathrm{R})$, where $\mathrm{R}$ is
the set of regular values of the energy momentum map.

\section{The classical spherical pendulum}
\label{sec2} 
%%%%%%%%%%%%%%

In this section we describe the geometry of the
classical spherical pendulum. More details can be found 
in \cite[chpt V]{cushman-bates}. \medskip 

\subsection{The basic system}
\label{sec2subsec1}
%%%%%%%%%%%%%%%%%%

We discuss the spherical pendulum as a constrained system. First we give the unconstrained system. Let 
$T^{\ast }{\mathbb{R} }^3 = {\mathbb{R}}^3 \times ({\mathbb{R}}^3)^{\ast }$ have coordinates $(q,p)$ 
and symplectic form $\widetilde{\omega } = \sum^3_{i=1} \dee p_i
\wedge \dee q_i = \dee \theta $, where $\theta = \langle p, \dee q \rangle $. Here $\langle \, \, , \, \, \rangle $ 
is the Euclidean inner on ${\R }^3$, which we use to identify $T^{\ast }{\R }^3$ with $T{\R }^3$. The unconstrained Hamiltonian system $(H, T{\R} ^3, \widetilde{\omega })$ has unconstrained Hamiltonian 
\begin{equation*}
\widetilde{H}:T{\R}^3 \rightarrow \R : (q,p) \mapsto %
\onehalf \langle p,p \rangle +\langle q, e_3 \rangle .
\end{equation*}
Here ${\{ e_i \} }^3_{i=1}$ is the standard basis of ${\R} ^3$. Now constrain the system 
$(\widetilde{H},T{\mathbb{R}}^{3},\widetilde{\omega })$ to the tangent bundle 
$TS^{2}=\{(q,p)\in T{\mathbb{R}}^{3}\,\mathop{\rule[-4pt]{.5pt}{13pt}\, }%
\nolimits\,\langle q,q\rangle =1\,\,\&\,\,\langle q,p\rangle =0\} $ of the $2$-sphere $S^{2}$ with symplectic form 
$\omega =\widetilde{\omega }_{\mid TS^{2}}$. Again we use the Riemannian metric on $S^2$ induced from 
the Euclidean inner product on ${\R }^3$ to identify the cotangent bundle $T^{\ast }S^2$ with the 
tangent bundle $TS^2$. The constrained Hamiltonian is $H=\widetilde{H}_{\mid TS^{2}}$, 
that is, 
\begin{equation}
H:TS^{2}\subseteq T{\R }^3 \rightarrow \mathbb{R}:(q,p)\mapsto 
\mbox{$\frac{\scriptstyle
1}{\scriptstyle 2}\,$}\langle p,p\rangle +\langle q,e_{3}\rangle .
\label{eq-s3ss1one}
\end{equation}%
The classical spherical pendulum is the Hamiltonian system 
$(H,TS^{2},\omega )$. The integral curves of the Hamiltonian vector field $X_H$ of the
Hamiltonian $H$ (\ref{eq-s3ss1one}) satisfy 
\begin{subequations}
\begin{align}
\frac{\dee q}{\dee t} & =p 
\label{eq-s3ss1twoa} \\ 
\rowspace \frac{\dee p}{\dee t} & =-e_{3}+(\langle q,e_{3}\rangle -\langle p,p\rangle )q%
\label{eq-s3ss1twob}
\end{align}
\end{subequations}
on $T{\R}^{3}$. Since $TS^{2}$ is an invariant manifold of 
(\ref{eq-s3ss1twoa})--(\ref{eq-s3ss1twob}), it follows that they define the integral curves of a vector field $X_H = X_{\widetilde{H}}|TS^{2}$, which governs the motion of the
spherical pendulum. A calculation shows that $H=\widetilde{H}_{\mid TS^{2}}$ and $L=%
\widetilde{L}_{\mid TS^{2}}$ are constants of motion of the vector field 
$X_H$. 

\subsection{Reduction of symmetry}
\label{sec2subsec2} 
%%%%%%%%%%%%%%

The angular momentum $\widetilde{L}$ of the unconstrained system $(%
\widetilde{H}, T{\mathbb{R} }^3, \widetilde{\omega })$ is a constant of
motion, because the unconstrained Hamiltonian $\widetilde{H}$ is invariant
under the $S^1$-action 
\begin{equation*}
\widetilde{\Phi }: S^1 \times T{\R}^3 \rightarrow T{\R}^3: %
\big( s, (q,p) \big) \mapsto {\widetilde{\Phi }}_s(q,p) = (R_s q, R_sp),
\end{equation*}
where $R_s=${\tiny $%
\begin{pmatrix}
\cos s & -\sin s & 0 \\ 
\sin s & \cos s & 0 \\
0 & 0 & 1
\end{pmatrix}%
$}. Since $R_s$ is a rotation, $TS^2$ is invariant under ${\widetilde{\Phi }}_s$. The 
infinitesimal generator of the $S^1$-action restricted to $TS^2$ is $X_{\widetilde{L}}|TS^2$, whose integral 
curves satisfy
\begin{subequations}
\begin{align}
\frac{\dee q}{\dee s} & = -q \times e_3
\label{eq-s3ss2onea} \\
\rule{0pt}{19pt} \frac{\dee p}{\dee s} & = - p \times e_3. 
\label{eq-s3ss2oneb}
\end{align}
\end{subequations}
So the constrained Hamiltonian $H$ is invariant under the $S^1$-action 
\begin{equation}
\Phi : S^1 \times TS^2 \rightarrow TS^2: \big( s, (q,p) \big) \mapsto {\Phi }%
_s(q,p) = (R_s q, R_sp).  
\label{eq-s3ss2one}
\end{equation}
The invariance of the constrained angular momentum $L = \widetilde{L}|TS^2$ under the $S^1$-symmetry ${\Phi }_s$ shows that it is an integral of the spherical pendulum. Thus the spherical pendulum is an integrable system 
$(H, L, TS^2, \omega )$. \medskip

Using invariant theory we reduce the $S^{1}$-symmetry of the spherical
pendulum. First observe that the algebra of polynomials on $T{\R}^{3}$, 
which are invariant under the $S^{1}$-action $\widetilde{\Phi }$, is
generated by 
\begin{subequations}
\begin{equation}
\begin{array}{rlcrlcrl}
{\pi }_{1} & =q_{3}, & \quad  & {\pi }_{2} & =p_{3} & \quad  & {\pi }_{3} & 
=p_{1}^{2}+p_{2}^{2}+p_{3}^{2} \\ 
\rule{0pt}{16pt}{\pi }_{4} & =q^2_1+q^2_2 & \quad  & {\pi }_{5} & 
=q_{1}p_1+q_{2}p_2 & \quad  & {\pi }_{6} & =q_{1}p_{2}-q_{2}p_{1},%
\end{array}
\label{eq-s3ss2twoa}
\end{equation}%
subject to the relation 
\begin{equation}
{\pi }_{5}^{2}+{\pi }_{6}^{2}={\pi }_{4}({\pi }_{3}-{\pi }_{2}^{2}),%
\mspace{10mu}\text{where}\mspace{5mu } ({\pi }_{3} - {\pi }_{2}^{2}) \geq 0%
\mspace{5mu}\text{and}\mspace{5mu} {\pi }_{4}\geq 0.  \label{eq-s3ss2twob}
\end{equation}%
The algebra of polynomials invariant under the $S^{1}$-action $\Phi $ 
(\ref{eq-s3ss2one}) on $TS^{2}$ is generated by $({\pi }_{i})|_{TS^{2}}$ for 
$1\leq i\leq 6$, which we will denote by ${\pi }_{i}$, subject to the
additional relations 
\begin{equation}
\begin{array}{rl}
{\pi }_{4}+{\pi }_{1}^{2} & =1 \\ 
\rule{0pt}{12pt}{\pi }_{5}+{\pi }_{1}{\pi }_{2} & =0 . 
\end{array}
\label{eq-s3ss2twoc}
\end{equation}%
These relations are just the defining equations of $TS^{2}$ expressed in terms of
invariants. Eliminating ${\pi }_{4}$ and ${\pi }_{5}$ from the relation (\ref%
{eq-s3ss2twob}) using (\ref{eq-s3ss2twoc}) gives 
\end{subequations}
\begin{equation}
{\pi }_{2}^{2}+{\pi }_{6}^{2}={\pi }_{3}(1-{\pi }_{1}^{2}),\mspace{10mu}%
\text{where}\mspace{5mu}-1\leq {\pi }_{1}\leq 1\mspace{5mu}\text{and}%
\mspace{5mu} {\pi }_{3}\geq 0,  
\label{eq-s3ss2three}
\end{equation}%
which defines $W=TS^{2}/S^{1}$, the space of $S^{1}$ orbits of the
action $\Phi $ on $TS^{2}$. In terms of invariants the space $L^{-1}({\ell })
$ is defined by ${\pi }_{6}=\ell $. Thus after removing the $S^{1}$ symmetry
of the spherical pendulum, the \emph{reduced phase space} $P_{\ell
}=L^{-1}(\ell )/S^{1}$ is the subvariety of ${\R}^{3}$ with
coordinates $({\pi }_{1},{\pi }_{2},{\pi }_{3})$ given by 
\begin{equation}
{\pi }_{2}^{2}+{\ell }^{2}={\pi }_{3}(1-{\pi }_{1}^{2}),\mspace{10mu}\text{%
where}\mspace{5mu}-1\leq {\pi }_{1}\leq 1\mspace{5mu}\text{and}\mspace{5mu}
{\pi }_{3}\geq 0.  
\label{eq-s3ss2four}
\end{equation}

Since the Hamiltonian $H$ of the spherical pendulum is invariant under the $%
S^1$-action $\Phi $, it induces the \emph{reduced Hamiltonian} 
\begin{equation}
{\widehat{H}}_{\ell}: P_{\ell } \subseteq {\R }^3 \rightarrow \R: ({\pi }%
_1, {\pi }_2, {\pi }_3) \mapsto \onehalf  {\pi }_3 + {\pi }_1.
\label{eq-s3.2fournew}
\end{equation}
When $\ell \ne 0$, then $-1< {\pi }_1 < 1$. So equation (\ref{eq-s3ss2four})
may be written as ${\pi }_3 = \frac{{\pi }^2_2 + {\ell }^2}{1-{\pi }^2_1}$, where $ |{\pi }_1| <1$.   
Thus $P_{\ell }, \, \, \ell \ne 0$, is diffeomorphic to ${\widehat{\R}}^2 = (-1,1)\times \R $ with coordinates 
$({\pi }_1, {\pi }_2)$ via the diffeomorphism 
\begin{equation}
\psi : {\widehat{\R}}^2 \subseteq {\R}^2 \rightarrow
P_{\ell } \subseteq {\R}^3:({\pi }_1, {\pi}_2) \rightarrow  ({\pi }_1, {\pi }_2, \frac{{\pi }^2_2 + 
{\ell }^2}{1-{\pi }^2_1}).  
\label{eq-s3ss2sixstar}
\end{equation}
On ${\widehat{\R}}^2$ the reduced Hamiltonian becomes 
\begin{equation}
H_{\ell } = {\psi }^{\ast }{\widehat{H}}_{\ell }: {\widehat{\R}}^2
\rightarrow \mathbb{R}: ({\pi }_1, {\pi }_2) \mapsto \ttfrac{1}{2}  
\ttfrac{1}{1-{\pi }^2_1} {\pi }^2_2 + V_{\ell }({\pi }_1),
\label{eq-s3ss2seven}
\end{equation}
where $V_{\ell }({\pi }_1) = {\pi }_1 + 
\mbox{$\frac{{\scriptstyle {\ell}^2}}{{\scriptstyle 2(1-{\pi }^2_1)}}$}$, $\ell \ne 0$. \medskip

To complete the reduction process, we determine the symplectic
structure on ${\widehat{R}}^2$ and the reduced equations of motion. We begin with finding the symplectic
form on the reduced phase space. The Poisson
bracket $\mbox{${\{ \, \,  , \, \, \}} _{{\scriptstyle {\R }^6}}$}$  has structure matrix 
${\mathcal{W}}_{{\R}^6}$  given in table 2.1. \bigskip   

\noindent \hspace{.25in}%
\begin{tabular}{c|cccccc|c}
$\mbox{${\{ A , B \}} _{{\scriptstyle {\R}^6}}$}$ & ${\pi }_1$ & ${\pi}_2$ & 
${\pi }_3$ & ${\pi }_4$ & ${\pi }_5 $ & ${\pi }_6$ & $B$ \\ \hline
${\pi }_1$ & $0$ & $1$ & $2{\pi }_2 $ & $0$ & $0$ & $0$ &  \\ 
${\pi }_2$ & $-1$ & $0$ & $0$ & $0$ & $0$ & $0$ &  \\ 
${\pi }_3$ &$-2{\pi}_2$  &$0$  & $0$ & $-4{\pi }_5$ & $2({\pi }^2_2-{\pi }_3)$ & $0$ &  \\ 
${\pi }_4$ & $0$ &$0$  & $4{\pi }_5$ & $0$ & $2{\pi }_4$ & $0$ &  \\ 
${\pi }_5$ &$0$  &$0$  & $-2({\pi }^2_2-{\pi }_3)$ &$-2{\pi }_4$ & $0$ & $0$ &  \\ 
${\pi }_6$ & $0$ &$0$  & $0$ & $0$ & $0$ & $0$ &  \\ \hline
$A$ &  &  &  &  &  &  & 
\end{tabular}
\bigskip

\noindent \hspace{1in}{Table 2.1. Structure matrix ${\mathcal{W}}_{{\R}^6}$ for 
$\mbox{${\{ \, \,  , \, \,  \}} _{{\scriptstyle {\R}^6}}$}$.} \bigskip 

Using table 2.1 we find that ${\{ {\pi }_4 +{\pi }^2_1 , {\pi }_5 +{\pi }_1 {\pi }_2 \} }_{\scriptstyle {\R}^6} |_{W}
 = 2$. So we may use Dirac brackets to compute the Poisson bracket ${\{ \, \,  , \, \,  \}}_{\scriptstyle W}$ 
 on the orbit space $W =TS^1/S^1$. 
We obtain the skew symmetric structure matrix ${\mathcal{W}}_{W}$ for the Poisson
bracket ${\{ \, \,  , \, \,  \}} _{\scriptstyle W}$ is given in table 2.2. \bigskip 

\noindent \hspace{1in}%
\begin{tabular}{c|cccc|c}
$\mbox{${\{ A , B \}} _{{\scriptstyle W}}$}$ & ${\pi }_1$ & ${\pi}_2$ & ${%
\pi }_3$ & ${\pi }_6$ & $B$ \\ \hline
${\pi }_1$ & $0$ & $1-{\pi }^2_1$ & $2{\pi }_2 $ & $0$ &  \\ 
${\pi }_2$ &$-(1-{\pi }^2_1)$  & $0$ & $-2{\pi }_1{\pi }_3$ & $0$ &  \\ 
${\pi }_3$ &$-2{\pi }_2$  & $2{\pi }_1{\pi }_3$ & $0$ & $0$ &  \\ 
${\pi }_6$ & $0$ &$0$  & $0$ & $0$ &  \\ \hline
$A$ &  &  &  &  & 
\end{tabular}
\bigskip

\noindent \hspace{1.25in}{Table 2.2. Structure matrix ${\mathcal{W}}_{W}$
for $\mbox{${\{ \, \,  , \, \,  \}} _{{\scriptstyle W}}$}$.} \bigskip

\noindent  Because ${\pi }_6$ Poisson commutes with every smooth function on $W$, the structure matrix ${\mathcal{W}}_{P_{\ell }}$ for the Poisson bracket $%
\mbox{${\{ \, \, , \, \,  \}} _{{\scriptstyle P_{\ell }}}$}$ on $P_{\ell }$ is given in 
table 2.3. \bigskip

\noindent \hspace{1in}%
\begin{tabular}{c|ccc|c}
$\mbox{${\{ A , B \}} _{{\scriptstyle P_{\ell }}}$}$ & ${\pi }_1$ & ${\pi}_2$
& ${\pi }_3$ & $B$ \\ \hline
${\pi }_1$ & $0$ & $1-{\pi }^2_1$ & $2{\pi }_2 $ &  \\ 
${\pi }_2$ &$-(1-{\pi }^2_1)$  & $0$ & $-2{\pi }_1{\pi }_3$ &  \\ 
${\pi }_3$ &$-2{\pi}_2$  &$2{\pi}_1{\pi}_3$  & $0$ &  \\ \hline
$A$ &  &  &  & 
\end{tabular}
\medskip 

\noindent \hspace{1.25in}{Table 2.3. Structure matrix ${\mathcal{W}}%
_{P_{\ell}}$ for $%
\mbox{${\{ \, \,  , \, \,  \}} _{{\scriptstyle P_{\ell}}}$}$.} \bigskip

Thus the reduced equations of motion on the reduced phase space $P_{\ell }$ are 
\begin{subequations}
\begin{align}
{\dot{\pi }}_1 & = \mbox{${\{ {\pi }_1 , {\widehat{H}}_{\ell }
\}} _{{\scriptstyle P_{\ell }}}$} = \onehalf  
\mbox{${\{ {\pi }_1 , {\pi }_3 \}} _{{\scriptstyle
P_{\ell }}}$} + \mbox{${\{ {\pi }_1 , {\pi }_1 \}} _{{\scriptstyle P_{\ell}}}$} 
= {\pi }_2  \label{eq-s3ss2eighta} \\
{\dot{\pi }}_2 & = 
\mbox{${\{ {\pi }_2 , {\widehat{H}}_{\ell } \}} _{{\scriptstyle P_{\ell }}}$} = 
\onehalf \mbox{${\{ {\pi }_2 , {\pi }_3 \}} _{{\scriptstyle P_{\ell }}}$} + 
\mbox{${\{ {\pi }_2 , {\pi }_1 \}} _{{\scriptstyle P_{\ell}}}$} 
= -{\pi }_1{\pi }_3 +{\pi }^2_1 -1  
\label{eq-s3ss2eightb} \\
{\dot{\pi }}_3 & = 
\mbox{${\{ {\pi }_3 , {\widehat{H}}_{\ell }
\}} _{{\scriptstyle P_{\ell }}}$} = \onehalf \mbox{${\{ {\pi }_3 , {\pi }_3 \}} _{{\scriptstyle
P_{\ell }}}$} + \mbox{${\{ {\pi }_3 , {\pi }_1 \}} _{{\scriptstyle P_{\ell}}}$} = -2{\pi }_2.  
\label{eq-s3ss2eightc}
\end{align}
\end{subequations}
Note that the function ${\pi}^2_2 +{\ell }^2 - {\pi }_3(1-{\pi }^2_1)$, whose zero set 
defines the \linebreak 
reduced phase space $P_{\ell }$, is a constant of motion of 
the reduced equations of motion. On ${\widehat{R}}^2$ with coordinates 
$({\pi }_1, {\pi }_2) $ the structure matrix ${\mathcal{W}}_{{\widehat{R}}^2}$ 
of the Poisson bracket $\mbox{${\{ \, \,  , \, \,  \}} _{{\scriptstyle {\widehat{R}}^2}}$}
= {\psi }^{\ast }\mbox{${\{ \, \,  , \, \,  \}} _{{\scriptstyle P_{\ell}}}$}$
is {\tiny $%
\begin{pmatrix}
0 & 1-{\pi }^2_1 \\ 
-(1-{\pi }^2_1) & 0%
\end{pmatrix}
$}, which is invertible on $P_{\ell }$, $\ell \ne 0$ since $1-{\pi }^2_1 > 0$. Thus  
the symplectic form on ${\widehat{\R} }^2$ is 
${\omega }_{{\widehat{\R}}^2} = \mbox{$\ttfrac{{\scriptstyle
1}}{{\scriptstyle 1-{\pi }^2_1}}$} \, \dee {\pi }_2 \wedge \dee {\pi }_1$, which 
corresponds to the matrix $({\mathcal{W}}^{-1}_{{\widehat{R}}^2})^T=${\tiny 
$\frac{1}{1-{\pi }^2_1} \begin{pmatrix} 0 & 1 \\ -1 & 0 \end{pmatrix}$}. 

So the reduced equations of motion on ${\widehat{\R}}^2$ of the reduced
system $(H_{\ell }, {\widehat{\R}}^2, $ ${\omega }_{{\widehat{\R}}^2})$ with $\ell \ne 0$ are 
\begin{subequations}
\begin{align}
{\dot{\pi }}_1 & = 
\mbox{${\{ {\pi }_1 , H_{\ell } \}} _{{\scriptstyle{\widehat{\R}}^2}}$} = 
\mbox{${\{ {\pi }_1 , \onehalf 
\ttfrac{1}{1-{\pi}^2_1}\, {\pi }^2_2 \}} _{{\scriptstyle {\widehat{\R }}^2}}$} + 
\mbox{${\{ {\pi }_1 , V_{\ell }({\pi }_1) \}} _{{\scriptstyle {\widehat{\R }}^2}}$}  \notag \\
& = \mbox{$\frac{{\scriptstyle {\pi }_2}}{{\scriptstyle 1-{\pi }^2_1}}$}\, 
\mbox{${\{ {\pi }_1 , {\pi }_2 \}} _{{\scriptstyle {\widehat{\R }}^2}}$} = 
{\pi }_2  
\label{eq-s3ss2ninea} \\
{\dot{\pi }}_2 & = \mbox{${\{ {\pi }_2 , H_{\ell } \}} _{{\scriptstyle{\widehat{R}}^2}}$} = 
\mbox{${\{ {\pi }_2 , \onehalf \ttfrac{1}{1-{\pi}^2_1}\, {\pi }^2_2 \}} _{{\scriptstyle {\widehat{\R }}^2}}$} + \mbox{${\{{\pi }_2 , V_{\ell }({\pi }_1) \}} _{{\scriptstyle {\widehat{\R }}^2}}$} \notag \\
& = \mbox{$\frac{{\scriptstyle {\pi }_1{\pi }^2_2}}{{\scriptstyle (1-{\pi}^2_1)^2}}$}\, 
\mbox{${\{ {\pi }_2 , {\pi }_1 \}} _{{\scriptstyle{\widehat{\R}}^2}}$} 
+V^{\prime}_{\ell }({\pi}_1)\poisbr{{\pi }_2}{{\pi }_1}{{\widehat{\R}}^2} \notag \\
& = \mbox{$-\frac{{\scriptstyle {\pi }_1{\pi}^2_2}}{{\scriptstyle 1-{\pi }^2_1}}$} -
V^{\prime }_{\ell }({\pi }_1)(1-{\pi }^2_1).  
\label{eq-s3ss2nineb} 
\end{align}
\end{subequations}
Note that the reduced Hamiltonian $H_{\ell }$ is a constant of motion
of the reduced equations of motion.

\subsection{Regular values}
\label{sec2subsec3} 
%%%%%%%%%%%%%

In this section we determine the topology of the set of regular values,
which lie in the image of the energy-momentum mapping 
\begin{equation*}
\mathcal{EM}:TS^2 \rightarrow {\mathbb{R} }^2:(q,p) \mapsto \big( H(q,p),
L(q,p) \big) = \big( \onehalf \langle p,p \rangle +\langle q, e_3 \rangle , 
\, q_1p_2-q_2p_1\big)
\end{equation*}
of the spherical pendulum. \medskip

First we determine the set of critical values of the energy momentum map.
The pair $(h,\ell )$ is a critical value of $\mathcal{EM}$ if and only if the 
$h$-level set of the reduced Hamiltonian ${\widehat{H}}_{\ell }$, that
is, the $2$-plane $\onehalf {\pi }_3 +{\pi }_1 = h$ in ${\mathbb{R} }^3$ with
coordinates $({\pi }_1, {\pi }_2, {\pi }_3)$, intersects the reduced space 
$P_{\ell } \subseteq {\R}^3$, defined by ${\pi }^2_2 +{\ell }^2 = 
{\pi }_3(1-{\pi }^2_1)$ with $|{\pi }_1| \le 1$ and ${\pi }_3 \ge 0$, at a
point of multiplicity greater than $1$. In other words, the polynomial 
\begin{displaymath}
Q({\pi }_1, {\pi }_2) = {\pi }^2_2-\big( 2(h-{\pi }_1)(1-{\pi }^2_1) - {\ell }^2 \big) = {\pi }^2_2 - P_{h,\ell }({\pi }_1),
\end{displaymath}
which is obtained by eliminating ${\pi }_3$ from the defining equation of $P_{\ell }$, has a
multiple root $({\pi }_1, {\pi }_2) \in [-1,1] \times \R $, that is, $0= Q({\pi }_1, {\pi }_2) $ and 
$(0,0) = DQ({\pi }_1, {\pi }_2) = \big( -P^{\prime }_{h, \ell }({\pi }_1), 2{\pi }_2 \big)$. Clearly ${\pi }_2 =0$ 
and ${\pi }_1$ is a multiple root of $P_{h, \ell }$ in $[-1,1] $. Let $\Delta $ be the discriminant of $P_{h,\ell }$, that is, 
$\Delta = \gcd (P_{h, \ell }, P^{\prime }_{h,\ell })$. The set of all $ (h, \ell ) \in {\mathbb{R} }^2$ such that $P_{h, \ell }$ has a multiple root in $[-1,1]$ is the discriminant locus $\{ \Delta =0 \} $ of $P_{h, \ell }$.
Suppose that $(h, \ell ) \in \{\Delta =0 \}$, then for some $s\in [-1,1]$ and $t\in \R $ we may write 
$P_{h, \ell }({\pi }_1) = 2({\pi }_1-s)^2({\pi }_1 -t)$.
Equating the coefficients of like powers of ${\pi }_1$ in the preceding equality gives 
\begin{align}
h & = 2s +t  \notag \\
-1 & = s^2+2st  
\label{eq-s3ss3one} \\
2h - {\ell }^2 & = -2ts^2. 
\notag
\end{align}
Eliminating $t$ from (\ref{eq-s3ss3one}) gives the following parametrization
of the discriminant locus $\{ \Delta =0 \} $ 
\begin{equation*}
(h(s), \ell (s) ) = \big( \ttfrac{3}{2} s - \ttfrac{1}{2}\ttfrac{1}{s}, 
\pm  (1-s^2)   \ttfrac{1}{\sqrt{-s}}  \big) , \mspace{10mu }\text{where} %
\mspace{5mu} s \in [-1, 0) \cup \{ 1 \} ,  
\end{equation*} 
see figure 1. Thus $\{ \Delta =0 \} $ is the union of two curves ${\mathcal{B}}_{\pm}$, 
which go to $+\infty $ as $s \nearrow 0$, are reflections in the $h$-axis of each
other, meet at a right angle and end when $s=-1$, that is, when $(h, \ell )
=(-1, 0)$. Otherwise they do not intersect. When $s=1$ we obtain the \emph{isolated point} 
$(1, 0)$. The image of the energy momentum mapping $\mathcal{EM}$
of the spherical pendulum is the closed subset of ${\R}^2$ bounded
by the curves ${\mathcal{B}}_{\pm }$, which contains the point $(1,0)$. 
The set $\mathrm{R}$ of regular values in the
image of $\mathcal{EM}$ is the interior of the image of $\mathcal{EM}$ with
the point $(1,0)$ removed. Thus $\mathrm{R}$ is diffeomorphic to an open $2$-disk
with its center deleted and so is not simply connected.  
\bigspace
\par \noindent \hspace{1in}\begin{tabular}{l}
\setlength{\unitlength}{2pt}
\vspace{-.9in}
\includegraphics[width=200pt]{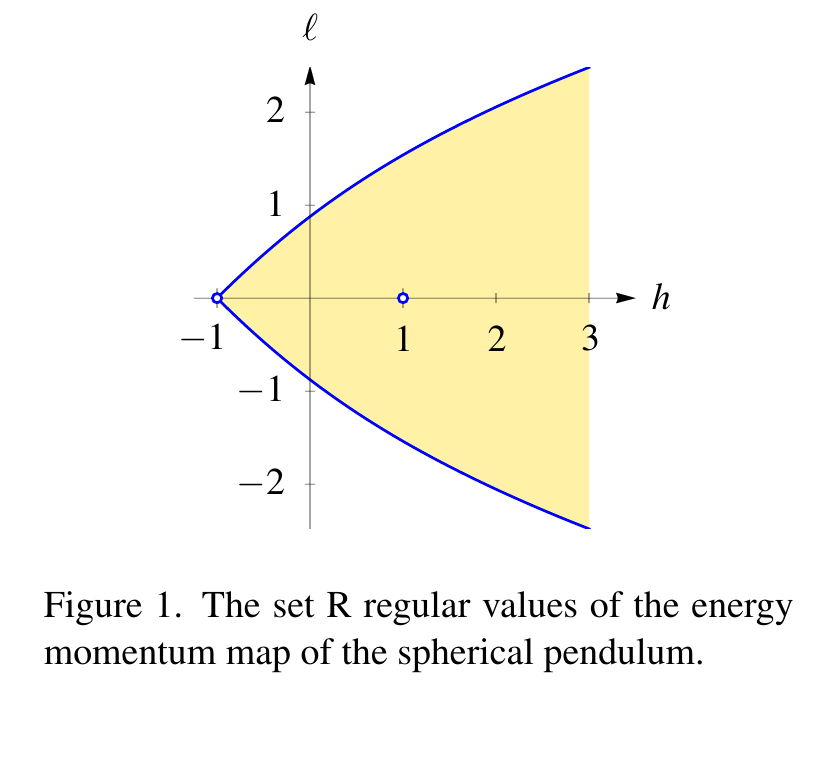} 
\vspace{.75in}
\end{tabular}

Actually, the argument above gives the following statification of the range of the energy 
momentum map $\mathcal{EM}$, which is a semialgebraic subset of ${\R }^2$. \medskip 

1. \parbox[t]{4.25in}{The set of regular values $\mathrm{R}$ is the top $2$-dimensional 
stratum. Here ${\mathcal{EM}}^{-1}(r)$ for $r \in \mathrm{R}$ is a smooth $2$-dimensional torus.} 
\smallspace 
\indent 2. \parbox[t]{4.25in}{The two curves ${\mathcal{B}}_{\pm }$, excluding the 
point $(-1,0)$ are two $1$-dimesional strata. Here ${\mathcal{EM}}^{-1}(b)$ for 
$b \in {\mathcal{B}}_{\pm}\setminus \{ (-1,0) \}$ is a smooth circle.} 
\vspace{.1in} 
\par 3. \parbox[t]{4.25in}{There are two $0$-dimensional strata: the points 
$(-1, 0)$ and $(1,0)$. Here ${\mathcal{EM}}^{-1}(-1,0)$ is a point; while ${\mathcal{EM}}^{-1}(1,0)$ is 
an immersed $2$-sphere in ${\R }^4$ with one normal crossing, in other words, a once 
pinched $2$-torus.} \medskip 

The orbit map of the $S^1$-action $\Phi $ on $L^{-1}(\ell )$ is given by 
\begin{equation}
\widetilde{\rho}:L^{-1}(\ell ) \subseteq TS^2 \rightarrow P_{\ell }\subseteq {\R} ^3: 
(q,p) \mapsto \big( {\pi }_1(q,p), \, {\pi }_2(q, p), \, {\pi }_3(q, p) \big) .
\label{eq-zeronew}
\end{equation}

Suppose that $(h, \ell ) \in \mathrm{R}$. Then the $h$-level set of the reduced Hamiltonian
 ${\widehat{H}}_{\ell}:P_{\ell } \subseteq {\R}^3 \rightarrow \R : 
({\pi }_1, {\pi }_2, {\pi }_3) \mapsto \onehalf {\pi }_3 + {\pi }_1$ is 
\begin{displaymath}
\mbox{\small $\big\{ \big({\pi }_1, {\pi }_2, 2(h-{\pi }_1) \big) \in P_{\ell } \, 
\mathop{\rule[-4pt]{.5pt}{13pt}\, }\nolimits \, {\pi }^2_2 = P_{h, \ell }({\pi }_1), \, \, \, 
\mbox{where $ -1 < {\pi }^{-}_1 \le {\pi }_1 \le  {\pi }^{+}_1 < 1$} \big\} $.}  
\end{displaymath}
Here ${\pi }^{\pm}_1$ are consecutive roots of the polynomial $P_{h, \ell }$ 
in $[-1,1]$. Thus $({\widehat{H}}_{\ell })^{-1}(h) = 
\widetilde{\rho }\big( H^{-1}(h) \cap L^{-1}(\ell ) \big)$ is a smooth submanifold of $P_{\ell }$, which is diffeomorphic to a circle. From the construction of the $S^1$-orbit space $P_{\ell }$ it follows that ${\widetilde{\rho}}^{-1}
\big( ({\widehat{H}}_{\ell })^{-1}(h) \big) $ is the total space of an $S^1$-bundle 
$\Sigma $ with base $S^1 = ({\widehat{H}}_{\ell })^{-1}(h)$. Because $P_{\ell }$ is
homeomorphic to ${\R} ^2$, we deduce that $({\widehat{H}}_{\ell })^{-1}(h)$ 
is contractible to a point in $P_{\ell }$. Thus $\Sigma $ is a
product bundle, which implies that ${\widetilde{\rho}}^{-1}
\big( ({\widehat{H}}_{\ell })^{-1}(h) \big) $ is diffeomorphic to a smooth $2$-torus 
$T^2_{h, \ell } = H^{-1}(h) \cap L^{-1}(\ell ) = {\mathcal{EM}}^{-1}(h, \ell )$. 
Consequently, we have
a smooth fibration 
\begin{equation}
\pi = \mathcal{EM}|_{{\mathcal{EM}}^{-1}(\mathrm{R})}: {\mathcal{EM}}^{-1}(\mathrm{R})\rightarrow \mathrm{R}
\label{eq-fifteennew}
\end{equation} 
whose fiber over $(h,\ell ) \in \mathrm{R}$ is the smooth torus $T^2_{h, \ell }$. Since every fiber of this 
mapping is compact, the fibration is locally trivial. Moreover, the bundle 
${\pi }|_{{\mathcal{B}}_{\pm}}: {\mathcal{EM}}^{-1}({\mathcal{B}}_{\pm }) \rightarrow {\mathcal{B}}_{\pm}$ 
is a trivial $S^1$ bundle over ${\R }$ because each 
curve ${\mathcal{B}}_{\pm}$ is contactible to a point. An integral curve of $X_L$ on 
${\mathcal{EM}}^{-1}({\mathcal{B}}_{\pm })$, which parametrizes a fiber of the bundle 
${\pi }|_{{\mathcal{B}}_{\pm}}$ is noncontractible curve and is 
the limit of an integral curve of $X_L$ on $T^2_{h,\ell }$ when $(h, \ell ) \in \mathrm{R}$ 
converges to a point on ${\mathcal{B}}_{\pm }$. 

\subsection{Action-angle coordinates}
\label{sec2subsec4} 
%%%%%%%%%%%%%%%%

In this section we construct local action-angle coordinates for the
spherical pendulum. \medskip

The smooth locally trivial $2$-torus fibration 
$\pi :{\mathcal{EM}}^{-1}(\mathrm{R})\rightarrow \mathrm{R}$ satisfies the hypotheses of the action-angle coordinate theorem for the integrable Hamiltonian system 
$(H,L,TS^{2},\omega )$ describing the spherical pendulum, see \cite[chpt IX]{cushman-bates}. 
Thus about each $(h, \ell )\in \mathrm{R}$ there is an open neighborhood $\mathcal{V}$ in $\mathrm{R}$ and 
an open neighborhood $\mathcal{U} = {\mathcal{EM}}^{-1}(\mathcal{V})$ in 
$TS^2$ of the $2$-torus $T_{h, \ell }^{2}={\mathcal{EM}}^{-1}(h,\ell)$ and a diffeomorphism 
\begin{equation*}
\varphi :\mathcal{U}\subseteq TS^{2}\rightarrow \varphi (\mathcal{U}) = 
\mathcal{V} \times {\mathbb{T}}^{2} \subseteq {\R }^2 \times {\mathbb{T}}^2 : (q,p)\mapsto 
(A_1, A_2, {\varphi }_{1},{\varphi }_{2}),
\end{equation*}%
where ${\mathbb{T}}^{2}$ is the affine $2$-torus ${\R}^{2}/(2\pi \Z )^{2}$, such
that 
\begin{itemize}
\item[1.] The symplectic form $\omega $ on $TS^2$ when restricted to 
$\mathcal{U}$ is exact and ${\varphi}_{\ast }({\omega }|_{\mathcal{U}}) = 
\sum^2_{i=1} \dee A_i \wedge \dee {\varphi }_i $.
\item[2.] The actions $A_1$ and $A_2$ are smooth functions of $H$ and $L$ on $\mathcal{U}$.
\item[3.] The vector field ${\varphi }_{\ast }(X_H|_{\mathcal{U}})$ is 
Hamiltonian on $(\mathcal{V} \times {\mathbb{T}}^2 ,
\sum^2_{i=1} \dee A_1 \wedge \dee {\varphi }_i )$, corresponding to the Hamiltonian 
$\mathcal{H} = {\varphi }_{\ast}(H|\mathcal{U})$. Moreover, the integral curves of $X_{\mathcal{H}}$ satisfy
\begin{align}
\frac{\dee {\varphi }_1}{\dee t} & = \frac{\partial \mathcal{H}}{\partial A_1} \hspace{.75in} \frac{\dee {\varphi }_2}{\dee t} = \frac{\partial \mathcal{H}}{\partial A_2}  
\notag \\
\largerowspace \frac{\dee A_1}{\dee t} & = -\frac{\partial \mathcal{H}}{\partial {\varphi}_1} = 0 \hspace{.4in} \frac{\dee A_2}{\dee t} = -\frac{\partial \mathcal{H}}{\partial {\varphi}_2} = 0 . \notag 
\end{align}
\end{itemize}
Choose $\mathcal{V}$ sufficiently small so that the fibration 
${\mathcal{EM}}|_{\mathcal{U}}$ is trivial. Then we have the following commutative 
diagram  
\begin{displaymath}
\begin{CD}
\mathcal{U} = {\mathcal{EM}}^{-1}(\mathcal{V}) @>\varphi>> \mathcal{V} \times {\mathbb{T}}^2  \\
@VV\mathcal{EM}|\mathcal{U}V                  @VV{\pi}V\\
\mathcal{V} \subseteq R @>(A_1,A_2)>> \mathcal{V} \subseteq {\R }^2
\end{CD}
\end{displaymath}
\medskip 

\noindent where $\pi : \mathcal{V} \times {\mathbb{T}}^2  \rightarrow \mathcal{V} \subseteq {\R }^2:
(A_1, A_2, {\varphi }_1, {\varphi }_2) \mapsto (A_1, A_2)$. Thus the bundle $\pi $ is a 
local trivialization of the bundle $\mathcal{EM}|_{{\mathcal{EM}}^{-1}(\mathrm{R})}$. \medskip

We now construct the action functions for the spherical pendulum. \linebreak 
According to the proof of the action-angle coordinate theorem in \cite{cushman-bates}, the action functions $A_{i}$ on 
$\mathcal{U}$ are constructed by finding linear combinations on $\mathcal{U}$ of the vector fields 
$X_{H|_{\mathcal{U}}}$ and $X_{L|_{\mathcal{U}}}$, whose coefficients are smooth functions on 
$\mathcal{EM}(\mathcal{U})=\mathcal{V}$, which have periodic flow of period 
$2\pi $ when restricted to $T_{h, \ell }^{2}$. Here $(h,\ell )=\mathcal{EM}(u)$
with $u\in \mathcal{U}$. Since $X_{L}$ has periodic flow of period $2\pi $
on ${\mathcal{EM}}^{-1}(\mathrm{R}) \subseteq TS^{2}$, we may define the second action 
function on $\mathcal{U}$ as $A_{2}= L$.  \medskip 

To construct the first action function $A_1$ we need to determine certain
functions related to the flows $({\varphi }^H_{t})|_{T_{h, \ell }^{2}}$ and 
$({\varphi }^L_{s})|_{T^2_{h, \ell }}$ of the vector fields $X_{H}$ and 
$X_{L}$ on $T_{h, \ell }^{2}$. For $u\in T_{h, \ell }^{2}$ let 
${\widetilde{\Gamma }}_{1}:[0,2\pi ]\rightarrow T_{h,\ell }^{2}:s\mapsto 
{\varphi }_{s}^{L}(u)$. Because the flow of $X_{L}$ on $T_{h, \ell }^{2}$ is periodic of
period $2\pi $, ${\widetilde{\Gamma }}_{1}$ is a closed
curve on $T_{h, \ell }^{2}$. Since $X_{H}(u)$ is nonzero and is transverse to 
$X_{L}(u)$ for every $u\in T_{h, \ell }^{2}$ and because $(h,\ell )$ is a
regular value of $\mathcal{EM}$, the integral curve ${\Gamma }_{2}:
t\mapsto {\varphi }_{t}^{H}(u)$, which starts at $u\in 
{\widetilde{\Gamma }}_{1}([0,2\pi ])$, has a unique positive first time $T$ for which 
${\varphi }_{T}^{H}(u)\in {\widetilde{\Gamma }}_{1}([0,2\pi ])$. The time $T$ is a
smooth function of $(h,\ell )$ and does not depend on the choice of starting
point $u$ in ${\widetilde{\Gamma }}_{1}([0,2\pi ])$, because it is the 
period of the reduced vector field $X_{{\widehat{H}}_{\ell }}$ on 
$({\widehat{H}}_{\ell })^{-1}(h)\subseteq P_{\ell }$. Let 
$\Theta $ be the smallest positive time it takes the curve $\widetilde{\Gamma }:s\mapsto 
{\varphi }_{-s}^{L}({\varphi }_{T}^{H}(u))$ to return to $u$, that is, 
${\varphi }_{-\Theta }^{L}({\varphi }_{T}^{H}(u))=u$. The function 
$\Theta $ does not depend on the choice of the point $u$ and is a smooth function of 
$(h, \ell )$. \medskip 

To find an explicit expression for the push forward ${\mathcal{A}}_1$ of first action $A_{1}$ to 
the image of the energy momentum mapping, consider the formula 
\begin{align}
{\mathcal{A}}_1(h, \ell ) & = \ttfrac{1}{2\pi } \, \int_{\Gamma } \langle p, \dee q \rangle |_{TS^2} 
=  \ttfrac{1}{2\pi} \, \int_{{\Gamma }_1} \langle p, \dee q \rangle |_{TS^2} +
 \ttfrac{1}{2\pi }\, \int_{{\Gamma }_2} \langle p, \dee q \rangle |_{TS^2},  \notag 
\end{align}
where $\Gamma : [0, 2\pi ] \rightarrow TS^2$ is a closed path in $T^2_{h,\ell }$ starting at 
$(q,p)$, which is a sum of two paths ${\Gamma }_1$ and ${\Gamma }_2$ on $T^2_{h, \ell }$.  
Suppose that $(h, \ell ) \in \mathrm{R} $. Choose the paths 
${\Gamma }_1$ and ${\Gamma }_2$ to be   
\begin{displaymath} 
{\Gamma }_1:[0,2\pi ] \rightarrow T^2_{h, \ell } \subseteq TS^2:
t' \mapsto {\varphi }^{\widetilde{T}(h, \ell )H}_{t'}(q,p) = \big( q(t'), p(t') \big), 
\end{displaymath}
and 
\begin{displaymath}
{\Gamma }_2:[0,2 \pi ] \rightarrow T^2_{h, \ell } \subseteq TS^2:s' \mapsto 
{\varphi }^{\widetilde{\Theta }(h, \ell ) L}_{-s'}(q,p) = \big( q(s'), p(s') \big) , 
\end{displaymath} 
respectively. Here $\widetilde{T}(h, \ell ) = T(h, \ell )/(2\pi )$ and 
$\widetilde{\Theta }(h, \ell ) = \Theta (h, \ell )/(2\pi )$. Using the path parameters as integration 
variables, we get 
\begin{align}
2\pi \, {\mathcal{A}}_1(h, \ell ) & = \int^{2\pi }_0 \langle p, \frac{\dee q}{\dee t' } \rangle \, \dee t' + 
\int^{2\pi }_0 \langle p, \frac{\dee q}{\dee s'} \rangle \, \dee s' \notag \\
& =  \int^{2\pi }_0 \langle p, \frac{\dee q}{\dee \, (t' \, \widetilde{T})} \rangle \, 
\dee \, (t' \, \widetilde{T})+ 
\int^{2\pi }_0 \langle p, \frac{\dee q}{\dee \, (-s' \widetilde{\Theta })} \rangle 
\, \dee \, (-s' \widetilde{\Theta }) \notag \\
& =  \int^{T(h, \ell )}_0 \langle p, \frac{\dee q}{\dee t} \rangle \, \dee t + 
\int^{-\Theta (h, \ell )}_0 \langle p, \frac{\dee q}{\dee s} \rangle \, \dee s ,  \notag 
\end{align} 
changing integration variable to the dynamical 
times $t = t' \, \widetilde{T}(h, \ell )$ and $s = -s' \widetilde{\Theta } (h, \ell )$. Here 
\begin{displaymath}
\frac{\dee }{\dee t}\hspace{-1pt} \begin{pmatrix} 
q \\ p \end{pmatrix} = X_H(q,p) \quad \mathrm{and} \quad  
\frac{\dee }{\dee s} \hspace{-1pt}\begin{pmatrix} q \\ p \end{pmatrix} = X_L(q,p). 
\end{displaymath}
Using (\ref{eq-s3ss1twoa}) and (\ref{eq-s3ss2onea}) we obtain 
\begin{align}
2\pi {\mathcal{A}}_1(h, \ell ) & = \int^{T(h, \ell )}_0 \langle p, p \rangle \, \dee t + 
\int^{-\Theta (h, \ell )}_0 \langle p, -q \times e_3 \rangle \, \dee s  \notag \\
& = \int^{T(h, \ell )}_0 {\pi }_3 \, \dee t  -\Theta (h, \ell ) \ell , \, \, \, 
\mbox{since $L = \langle q \times p , e_3 \rangle = \ell $ on $T^2_{h, \ell }$.} \notag 
\label{eq-s3ss4fournw}%
\end{align}
But $h = H_{\ell } = \onehalf {\pi }_3 + {\pi }_1$ and $\dee t = \frac{1}{{\pi }_2} \dee {\pi}_1 $, 
using (\ref{eq-s3ss2eighta}). So
\begin{align}
\int^{T(h, \ell )}_0 {\pi }_3 \, \dee t  & = 2\int^{T(h, \ell )}_0 (h-{\pi }_1) \, \dee t  \notag \\
& = 2h\, T(h, \ell )  - 4\int^{{\pi }^{+}_1}_{{\pi }^{-}_1} 
\frac{{\pi }_1}{\sqrt{2(h-{\pi }_1)(1-{\pi }^2_1) - {\ell }^2}} \dee {\pi }_1.  \notag  
\end{align}
Therefore, when $(h, \ell ) \in \mathrm{R}$, 
\begin{equation}
\mbox{\small $2\pi \, A_1(h, \ell ) = 2h\, T(h, \ell ) - 2\pi \, I(h, \ell) - \ell  \Theta (h, \ell ) $},  
\label{eq-s2ss4fiveavnw}
\end{equation}
where 
\begin{equation}
2\pi \, I(h, \ell ) = 4\int^{{\pi }^{+}_1}_{{\pi }^{-}_1} \frac{{\pi }_1}{\sqrt{2(h-{\pi }_1)(1-{\pi }^2_1) - {\ell }^2}} 
\dee {\pi }_1
\label{eq-s2ss4fsixvnw}
\end{equation}
Since 
\begin{subequations}
\begin{equation}
2\pi \widetilde{T}(h, \ell ) = T(h, \ell ) = 2\,  \int^{{\pi }^{+}_1}_{{\pi }^{-}_1}
 \frac{1}{\sqrt{2(h-{\pi }_1)(1-{\pi }^2_1) - {\ell }^2}} \dee {\pi }_1 
\label{eq-twentyfournwa}
\end{equation}
and
\begin{align}
2\pi \widetilde{\Theta }(h, \ell ) & = \Theta (h, \ell ) \notag \\
& = 2 \,  \ell \,  \int^{{\pi }^{+}_1}_{{\pi }^{-}_1}  
\frac{1}{(1-{\pi }^2_1)\sqrt{2(h-{\pi }_1)(1-{\pi }^2_1) - {\ell }^2}} \dee {\pi }_1, 
\label{eq-twentyfournwb}
\end{align}
\end{subequations}
we may rewrite the right hand side of (\ref{eq-s2ss4fiveavnw}) as 
\begin{align}
{\mathcal{A}}_1 (h, \ell ) & = 
\frac{1}{\pi } \int^{{\pi }^{+}_1}_{{\pi }^{-}_1} \frac{2(h-{\pi }_1)(1-{\pi }^2_1) - {\ell }^2}{(1-{\pi }^2_1)\sqrt{2(h-{\pi }_1)(1-{\pi }^2_1) - {\ell }^2}} \, \dee {\pi }_1 \notag \\
& = \frac{1}{\pi } \int^{{\pi }^{+}_1}_{{\pi }^{-}_1} \frac{\sqrt{2(h-{\pi }_1)(1-{\pi }^2_1) - {\ell }^2}}
{1-{\pi }^2_1} \, \dee {\pi }_1 .
\label{eq-s2ss4sixnw}
\end{align}
Using (\ref{eq-s2ss4sixnw}) straightforward calculation shows that 
\begin{equation}
\frac{\partial {\mathcal{A}}_1}{\partial h} = 
\widetilde{T}\, \, \, \mathrm{and} \, \, \,  \frac{\partial {\mathcal{A}}_1}{\partial \ell } = -\widetilde{\Theta }  
\label{eq-s2ss4sixstarnw}
\end{equation}
on $\mathrm{R} \setminus \{ \ell = 0 \} $. The reason for cutting the line $\{ \ell = 0 \} $ out of 
$\mathrm{R}$ is that the function $\widetilde{\Theta }$ has a jump discontinuity there, see fact 2.2. Consequently, 
$\frac{\partial {\mathcal{A}}_1}{\partial \ell }$ is not defined on $\mathrm{R} \cap \{ \ell = 0 \} $. \medskip 

\noindent \textbf{Proposition 2.1} On ${\mathcal{EM}}^{-1}(\mathrm{R})$ the function 
$A_1 = (\mathcal{EM})^{\ast }{\mathcal{A}}_1$ is an action. \medskip 

\noindent \textbf{Proof.} Suppose that $(h,\ell ) \in \mathrm{R} \setminus \{ \ell = 0 \} $. 
We have to show that the flow of the Hamiltonian vector field 
$X_{A_1}$ is periodic of period $2\pi $. From (\ref{eq-s2ss4sixstarnw}) it follows that 
$\dee \hspace{.5pt} {\mathcal{A}}_1 (h, \ell ) = 
\widetilde{T}(h,\ell ) \dee h - \widetilde{\Theta }(h, \ell ) \dee \ell $. Pulling back the preceding 
equation by the energy momentum mapping $\mathcal{EM}$ gives
\begin{displaymath}
\dee A_1 (q,p) = \widetilde{T}(h,\ell ) \dee H (q,p) - \widetilde{\Theta }(h, \ell ) \dee L (q,p), 
\end{displaymath}
for every $(q,p) \in T^2_{h, \ell } \subseteq TS^2$. Applying the map ${\omega }^{\flat}(q,p)$ yields  
\begin{equation}
X_{A_1}(q,p) = \widetilde{T}(h,\ell ) X_H(q,p) - \widetilde{\Theta }(h, \ell ) X_L (q,p)
\label{eq-sec2.5onevnw} 
\end{equation}
on $T^2_{h, \ell }$. On $T^2_{h, \ell}$ the flow ${\varphi }^H_t$ of $X_H$ commutes with the 
flow ${\varphi }^L_s$ of $X_L$. Thus the flow of $X_{A_1}$ on $T^2_{h, \ell }$ is given by 
\begin{displaymath}
{\varphi }^{A_1}_t (q,p) = 
{\varphi }^L_{-\widetilde{\Theta}(h,\ell )t} \comp {\varphi }^H_{\widetilde{T}(h, \ell )t} (q,p). 
\end{displaymath}
The preceding flow is periodic of period $2\pi $ because  
\begin{displaymath}
{\varphi }^{A_1}_{2\pi } (q,p)  = 
{\varphi }^L_{-2\pi \widetilde{\Theta}(h,\ell )} \comp {\varphi }^H_{2\pi \widetilde{T}(h, \ell )} (q,p)  
= {\varphi }^L_{-\Theta (h,\ell )} \comp {\varphi }^H_{T(h, \ell )} (q,p) = (q,p).  
\end{displaymath}
The last equality above follows from the definition of the time $T(h, \ell )$ of first return and the definition of the 
rotation number of the flow of $X_H$ on $T^2_{h,\ell }$. Therefore $A_1$ is an action function on 
$T^2_{h, \ell }$ with $(h, \ell ) \in \mathrm{R} \setminus \{ \ell = 0 \} $. \medskip 

For $(h,0) \in \mathrm{R} \cap \{ \ell =0 \} $ the action ${\mathcal{A}}_1$ (\ref{eq-s2ss4sixnw}) is well defined, 
but the Hamiltonian vector field $X_{A_1}$ on $T^2_{h,0}$ has a jump discontinuity along 
$T^2_{h,0} \cap L^{-1}(0)$ because for $(q,p) \in T^2_{h,0}$ 
\begin{displaymath}
X_{A_1}(q,p) = \left\{ \begin{array}{rl}
\widetilde{T}(h,0) X_H(q,p) \pm \onehalf X_L(q,p), & \mbox{if $-1 < h < 1$} \\
\smallrowspace \widetilde{T}(h,0) X_H(q,p) \pm  X_L(q,p), & \mbox{if $h > 1$, }
\end{array} \right. 
\end{displaymath}
see (\ref{eq-sec2.5onevnw}) and fact 2.2. However, the flow of $X_{A_1}$ on $T^2_{h,0}$ is 
\begin{displaymath}
{\varphi }^{A_1}_t = \left\{ \begin{array}{rl}
{\varphi }^H_{\widetilde{T}(h,0)t} \comp {\varphi }^L_{\onehalf t}, & \mbox{if $-1 < h < 1$} \\
{\varphi }^H_{\widetilde{T}(h,0)t} \comp {\varphi }^L_t, & \mbox{if $h > 1$.} 
\end{array} \right. 
\end{displaymath}
Indeed on $T^2_{h,0} \cap L^{-1}(0)$ we have 
\begin{displaymath}
{\varphi }^{A_1}_t = \left\{ \begin{array}{rl}
{\varphi }^H_{\widetilde{T}(h,0)t} \comp {\varphi }^L_{-\onehalf t}, & \mbox{if $-1 < h < 1$} \\
\smallrowspace {\varphi }^H_{\widetilde{T}(h,0)t} \comp {\varphi }^L_{-t}, & \mbox{if $h > 1$,} 
\end{array} \right. 
\end{displaymath}
because 
\begin{displaymath}
\left\{ \begin{array}{rl} 
{\varphi }^L_{-\onehalf t} = {\varphi }^{-L}_{-\onehalf t} = {\varphi }^L_{\onehalf t}, & \mbox{when $-1 < h < 1$} \\
{\varphi }^L_t = {\varphi }^{-L}_t = {\varphi }^L_t, & \mbox{when $h > 1$.} 
\end{array} \right. 
\end{displaymath}
Thus the flow ${\varphi }^{A_1}_t$ on $T^2_{h,0}$ is well defined and continuous for every 
$-1 < h < 1$ or $h > 1$. \medskip 

Now suppose that $(h,\ell ) \in \mathrm{R} \cap \{ \ell = 0 \} $. Since the bundle 
$\pi =\mathcal{EM}|_{{\mathcal{EM}}^{-1}(\mathrm{R}}: {\mathcal{EM}}^{-1}(\mathrm{R}) \rightarrow \mathrm{R}$ 
is locally trivial, for every $(h,0) \in \mathrm{R} \cap \{ \ell = 0 \}$ there is an open neighborhood $V$ such 
that ${\mathcal{EM}}^{-1}(V)$ is a trivial bundle. Thus there is a smooth section 
$s:V \rightarrow {\mathcal{EM}}^{-1}(V):(h, \ell ) \longmapsto s(h,\ell ) = \big( q(h, \ell ), p(h, \ell ) \big) $ 
such that $\pi (s(h,\ell )) = (h, \ell )$ for every $(h, \ell ) \in V$, that is, 
$\big( q(h, \ell ), p(h, \ell ) \big) \in T^2_{h, \ell }$. From the definition of the first action function we have 
\begin{equation}
{\varphi }^L_{2\pi \widetilde{\Theta }(h,\ell )} \comp 
{\varphi }^H_{2\pi \widetilde{T}(h, \ell )}\big( q(h, \ell ), p(h, \ell ) \big) = 
\big( q(h, \ell ), p(h, \ell ) \big) .
\label{eq-sec2.5twovnw}
\end{equation}
Taking the limit as $\pm \ell \nearrow 0$ gives
\begin{subequations}
\begin{equation}
{\varphi }^L_{\onehalf (2 \pi )} \comp {\varphi }^H_{T(h,0)} \big( q(h, 0 ), p(h, 0 ) \big)  = 
\big( q(h, 0 ), p(h, 0 ) \big) , 
\label{eq-sec2.5threeavnw} 
\end{equation}
if $-1 < h < 1 $ and 
\begin{equation}
{\varphi }^L_{2 \pi } \comp {\varphi }^H_{T(h,0)} \big( q(h, 0 ), p(h, 0 ) \big)  = 
\big( q(h, 0 ), p(h, 0 ) \big) , 
\label{eq-sec2.5threebvnw} 
\end{equation}
\end{subequations}
if $h > 1$. Consequently, the flow ${\varphi }^{A_1}_t$ is periodic of period $2\pi $ on $T^2_{h,0}$. 
This completes the proof that $A_1$ is an action function on ${\mathcal{EM}}^{-1}(\mathrm{R})$. 
\hfill $\square $ \medskip 

\subsection{Properties of the functions $\widetilde{\Theta }$, $\widetilde{T}$ and ${\mathcal{A}}_1$}
\label{sec2subsec5}
%%%%%%%%%%%%

We just state the principal analytic properties of the functions $\widetilde{\Theta }$ (\ref{eq-twentyfournwb}) and  
$\widetilde{T}$ (\ref{eq-twentyfournwa}). Their proofs may be found in 
\cite[chpt V and exercises]{cushman-bates}.  \medskip

\noindent \textbf{Fact 2.2} \medskip 

1. \parbox[t]{4.45in}{On $\mathrm{R} \setminus \{ \ell = 0 \}$ the function $\widetilde{\Theta }$ 
is real analytic and odd in $\ell $, that is, $\widetilde{\Theta }(h, - \ell ) = - \widetilde{\Theta }(h, \ell)$. 
Its principal value (\ref{eq-twentyfournwb}) has range $(-1, -\onehalf ) \cup (\onehalf , 1)$. } \smallspace
\indent 2. \parbox[t]{4.45in}{At $\mathrm{R} \cap \{ \ell = 0 \}$ the function $\widetilde{\Theta }$ has 
a jump discontinuity. Specifically, for $(h, \ell ) \in \mathrm{R} \cap \{ \ell > 0 \} $ 
\begin{equation}
\lim_{\ell \searrow 0}\widetilde{\Theta }(h, \ell ) = \left\{ 
\begin{array}{rl}
\onehalf , & \mbox{if $-1 < h < 1$} \\ 
1, & \mbox{if $h >1$.}
\end{array}
\right.  
\label{eq-s3ss5four}
\end{equation} }
\smallspace
\indent 3. \parbox[t]{4.45in}{The function $\widetilde{\Theta }:\mathrm{R} \subseteq {\R }^2 \rightarrow \R $ is  multivalued. Along any positively oriented closed curve, which generates the fundamental group of 
$\mathrm{R}$, its value decreases by $1$.} 
\smallspace
\indent 4. \parbox[t]{4.45in}{The function $\widetilde{T}:\mathrm{R} \subseteq {\R }^2 \rightarrow {\R }_{\ge 0}$ is real analytic and even in $\ell $, namely, $\widetilde{T}(h, - \ell ) = \widetilde{T}(h, \ell )$. 
Moreover, $\widetilde{T}(h,0) \nearrow \infty $ 
as $h \rightarrow 1$ and $\widetilde{T}(h,0) \searrow 0$ as $h \nearrow \infty$.}  \medskip 

We now determine the limiting values of the functions $\widetilde{T}$, $\widetilde{\Theta }$, and 
${\mathcal{A}}_1$ (\ref{eq-s2ss4sixnw}) as $(h, \ell ) \in \mathrm{R}$ converges to $(h(s), \ell (s))$ $\in 
\partial \overline{\mathrm{R}}$ for some $s \in [-1,0)$. \medskip 

\noindent \textbf{Proposition 2.3} Let $\widetilde{T}(s)= \widetilde{T}(h(s), \ell (s))$, 
$\widetilde{\Theta }(s) = \widetilde{\Theta }(h(s), \ell (s))$, and ${\mathcal{A}}_1(s) $ $= 
{\mathcal{A}}_1(h(s), \ell (s))$. Then 
\begin{subequations}
\begin{align}
& \widetilde{T}(s)  = \frac{\sqrt{-s}}{\sqrt{3s^2+1}}; 
\label{eq-s3ss4foura} \\
& \widetilde{\Theta }(s)  = \pm \frac{1}{\sqrt{3s^2+1}}, \, \, \, 
\mbox{when $\pm \ell (s) \ge 0$;} 
\label{eq-s3ss4fourb} \\
&  {\mathcal{A}}_1(s)  =0. 
\label{eq-s3ss4fourc}
\end{align}
\end{subequations} 

\noindent \textbf{Proof.} Consider ${\C }^{\vvee}$, 
the extended complex plane, which is cut along the real axis between 
$x_{-}$ and $x_{+}$ and again between $x_0$ and $\infty $. Here $x_{\pm ,0}$ are 
distinct roots of $P_{h, \ell }(x) = 2(h-x)(1-x^2) - {\ell }^2$ with 
\begin{displaymath}
\left\{ \begin{array}{rl}
-1 < x_{-} < x_{+} < h < 1 < x_0, & \mbox{if $-1 < h < 1$} \\
-1 < x_{-} < x_{+} < 1 < h < x_0, & \mbox{if $h > 1$}. \end{array} \right. 
\end{displaymath} 
Write $\sqrt{\onehalf P_{h, \ell }(z)} = 
\sqrt{r_{-}r_{+}r_0}\, {\mathrm{e}}^{i({\theta }_{-}+{\theta }_{+}+ 
{\theta }_0)/2}$, where $z-x_{0,\pm} = r_{0,\pm }{\mathrm{e}}^{i{\theta }_{0,\pm}}$ and 
$0 \le {\theta }_{0, \pm } < 2\pi $. For $(h, \ell ) = (h(s), \ell (s)) 
\in \partial \overline{\mathrm{R}}$ with 
$s \in [-1,0)$ the real polynomial $\onehalf P_{h,\ell }$ becomes 
\begin{displaymath}
\big( h(s) -x \big) \big( 1- x^2) - \onehalf {\ell }^2(s) = (x-s)^2(x-t), 
\end{displaymath}
where $t = t(s) = -\ttfrac{1}{s^2}(h(s) - {\ell }^2(s)) = -\onehalf s - \ttfrac{1}{2s}$. 
Let $\mathcal{C}$ be a positively 
oriented closed curve in ${\C }^{\vvee}$, which crosses the $\mathrm{Re}\, z$ 
axis twice: once in $(-1,x_{-})$ and once in $(x_{+}, 1)$. Note that the complex 
square root is negative just above the cut $[x_{-}, x_{+}]$. 
In the limit as $x_{\mp } \rightarrow s \in (-1,0)$ we see that $z - x_{\mp } \rightarrow 
z-x^{\ast }$, where $x^{\ast }= r^{\ast } \, {\mathrm{e}}^{i\, {\theta }^{\ast }}$ with 
$0 \le {\theta }^{\ast } < 2\pi $, because $r_{\mp} \rightarrow r^{\ast }$ and 
${\theta }_{\mp} \rightarrow {\theta }^{\ast }$. So 
$(z-s)\sqrt{z-t(s)} = 
\sqrt{r_0}\, {\mathrm{e}}^{\onehalf i\, {\theta }_0}(r^{\ast } \, {\mathrm{e}}^{i\, {\theta }^{\ast }})$. 
\medskip 

For $s \in [-1,0)$ we have 
\begin{align}
2 \pi \, \widetilde{\Theta } \big( h(s), \ell (s) \big) & = \frac{\ell (s)}{\sqrt{2}} \, 
\int_{\mathcal{C}} \frac{1}{(1-z^2)\sqrt{(h(s)-z)(1-z^2) - \ttfrac{1}{2}{\ell }^2(s)}} 
\, \dee z \notag \\
& = \frac{\ell (s)}{\sqrt{2}} \, \int_{\mathcal{C}} \frac{1}{(1-z^2)\sqrt{z-t}} \, 
\frac{\dee z}{z-s} \notag \\
& = \frac{\ell (s)}{\sqrt{2}} \Big( 2\pi i \, \frac{1}{(1-s^2)\sqrt{s-t}} \Big), \notag \\ 
&\hspace{.75in} \parbox[t]{3.25in}{because $\frac{1}{(1-z^2)\sqrt{z-t}}$ is complex analytic 
on ${\C}^{\vvee}$. Thus we can use Cauchy's integral formula.} \notag \\
& = \frac{2\pi i}{\sqrt{2}} \Big( \pm (1-s^2)\frac{1}{\sqrt{-s}} \Big) 
\frac{1}{(1-s^2)\sqrt{\ttfrac{3}{2}s + \ttfrac{1}{2s}}} = 
\pm \frac{2 \pi }{\sqrt{3s^2+1}}. \notag 
\end{align}   
Also for $s \in [-1,0)$ we have 
\begin{align}
2 \pi \widetilde{T}\big( h(s), \ell (s) \big) & = \frac{1}{\sqrt{2}}\, 
\int_{\mathcal{C}} \frac{1}{\sqrt{(h(s)-z)(1-z^2) - \ttfrac{1}{2}{\ell }^2(s)}} 
\, \dee z \notag \\
& = \frac{1}{\sqrt{2}} \, \int_{\mathcal{C}} \frac{1}{\sqrt{z-t}} \, \frac{\dee z}{z-s} 
= \frac{1}{\sqrt{2}} \Big( 2\pi i \frac{1}{\sqrt{s-t}} \Big) = 
\frac{2\pi \sqrt{-s}}{\sqrt{3s^2+1}}. \notag
\end{align}

\noindent We can write the integral $I(h, \ell ) = \frac{2}{\pi } \, 
\int^{x_{+}}_{x_{-}} \frac{x}{\sqrt{P_{h, \ell }(x)}} \, \dee x$  as the complex integral
\begin{equation}
I(h, \ell ) = \frac{1}{\pi \sqrt{2}} \int_{\mathcal{C}} \frac{z}{\sqrt{\onehalf P_{\ell ,h}(z)}} \, \dee z .
\label{eq-twodot}
\end{equation}
Thus when $(h,\ell )$ converges to $(h(s), \ell (s))$ the integral (\ref{eq-twodot}) becomes
\begin{align}
I(s) = I(h(s), \ell (s) ) & = \frac{1}{\pi \sqrt{2}} \int_{\mathcal{C}} \frac{z}{\sqrt{z-t(s)}} \, 
\frac{\dee z}{z-s} , 
\label{eq-three} \\
&\hspace{-.75in} = \frac{1}{\pi \sqrt{2}} \Big( 2\pi i \, \frac{s}{\sqrt{s-t(s)}} \Big) , \, \, \, 
\mbox{using Cauchy's integral theorem} 
\notag \\
& \hspace{-.75in} = 2i\, \frac{s}{i\sqrt{2}\sqrt{t(s)-s}}, \quad 
\mbox{by choice of complex square root} \notag \\
& \hspace{-.75in} =  -2\, \frac{(-s)^{3/2}}{\sqrt{3s^2+1}} . 
\label{eq-threedaggerstar}  
\end{align}
Therefore when $(h(s), \ell (s)) \in \partial \overline{\mathrm{R}} \cap \{ \ell \ge 0 \} $ 
using (\ref{eq-s2ss4fiveavnw}) we get 
\begin{align}
{\mathcal{A}}_1(s) & =  2h(s) \widetilde{T}(s) -\ell (s) \widetilde{\Theta }(s) - I(s) \notag \\
& = 2(\ttfrac{3}{2}s-\ttfrac{1}{2s})\frac{\sqrt{-s}}{\sqrt{3s^2+1}} 
-\frac{(1-s^2)}{\sqrt{-s}\, \sqrt{3s^2+1}} +2 \frac{(-s)^{3/2}}{\sqrt{3s^2+1}} \notag \\
& = \frac{1-3s^2 -(1-s^2) + 2s^2}{\sqrt{-s}\, \sqrt{3s^2+1}}  = 0.  
\label{eq-threestar}
\end{align} 
A similar argument shows that ${\mathcal{A}}_1(s) =0$ when $(h(s), \ell (s)) \in \partial \overline{\mathrm{R}} \cap 
\{ \ell \le 0 \} $. This confirms a result in \cite{richter-dullin-waalkens-wiersig}. \hfill $\square $  

\subsection{Additional properties of $A_1$}
\label{sec2subsec6}
%%%%%%%%%%%%

The next proposition gives some addtional properties of the first action function $A_1$ 
(\ref{eq-s2ss4sixnw}). \medskip 

\noindent \textbf{Proposition 2.4} The first action function ${\mathcal{A}}_1$ (\ref{eq-s2ss4sixnw}) is 
a positive real analytic function on $\mathrm{R} \setminus \{ \ell =0 \}$, which is 
continuous on $\mathrm{R}$, but is not differentiable on $\mathrm{R} \cap \{ \ell =0 \} $. Indeed, its partial derivative 
$\frac{\partial {\mathcal{A}}_1}{\partial \ell }$ has a jump discontinuity at 
$\mathrm{R}\cap \{ \ell = 0 \}$. Moreover, \medskip 

1. \parbox[t]{4.in}{${\mathcal{A}}_1$ has a nonnegative continuous extension to $\overline{\mathrm{R}}$, which vanishes on  $\partial \overline{\mathrm{R}}$.  Also ${\mathcal{A}}_1(-1,0) =0$, ${\mathcal{A}}_1(1,0) = 4/\pi $, 
and $\lim_{h\nearrow \infty} {\mathcal{A}}_1(h,0) = \infty $.}
\smallspace
\indent 2. \parbox[t]{4.in}{Let $\big( h(s_0), \pm \ell (s_0) \big)  \in 
\partial \overline{\mathrm{R}}$ for some $s_0 \in [-1,0)$. Then the function 
\begin{displaymath}
{\mathcal{A}}_1|_{\overline{\mathrm{R}} \cap \{ \ell = \pm \ell (s_0) \} }: [h(s_0), \infty ) \longmapsto 
\left[ 0, \infty ) \right. ,
\end{displaymath} 
is strictly increasing.} 
\smallspace
\indent 3. \parbox[t]{4.in}{Let $a > 0$. The $a$-level set of ${\mathcal{A}}_1$ in 
$\overline{\mathrm{R}}$ is the graph of a continuous, piecewise real analytic function 
\begin{displaymath}
{\widehat{A}}_a: \left[ -\infty, \infty ) \right. \rightarrow \left[ h^{\ast }(a), \infty ) \right. : \ell \longmapsto 
{\widehat{A}}_a(\ell ),
\end{displaymath} 
which is strictly decreasing when $\ell < 0$, is strictly increasing when $\ell > 0$, and has a positive minimum value 
$h^{\ast }(a)$ at $\ell =0$.} \medskip 

%\noindent See figure 2.
%\vspace{.15in}
%\bigspace
%\par \noindent \hspace{.75in}\begin{tabular}{l}
%\setlength{\unitlength}{2pt}
%\vspace{-.9in}
%\includegraphics[width=250pt]{levelsetsA1} 
%\vspace{.75in}
%\end{tabular}\medskip  

\noindent \textbf{Proof}. First we show that 
$I(h, \ell ) = \frac{2}{\pi }\int^{x_{+}}_{x_{-}} \frac{x}{\sqrt{2(h-x)(1-x^2) - {\ell }^2}} \, 
\dee x$ is a \linebreak 
locally real analytic function on $\mathrm{R}$. Write $I$ as the complex integral 
$I(h, \ell ) = \frac{1}{\pi } \int_{\mathcal{C}} \frac{z \, \dee z}{\sqrt{2(h-z)(1-z^2) - {\ell }^2}}$, where 
$\mathcal{C}$ and the complex square root are chosen as in the proof of (\ref{eq-s3ss4foura})--(\ref{eq-s3ss4fourc}).
Thinking of $h$ and $\ell $ as complex variables we obtain 
\begin{displaymath}
\frac{\partial I}{\partial \overline{h}} = \frac{1}{\pi } \int_{\mathcal{C}} \frac{\partial }{\partial \overline{h}} 
\left( \frac{z}{\sqrt{2(h-z)(1-z^2)-{\ell }^2}} \right) \dee z =0 
\end{displaymath}
and similarly $\frac{\partial I}{\partial \overline{\ell }} =0$. Therefore locally $I$ is a complex analytic function. 
Restricting $h$ and $\ell $ to be real variables shows that locally $I$ is a real analytic function on 
$\mathrm{R}$. Clearly, $I(h, -\ell ) = I(h, \ell )$. To show that $I$ is single-valued, consider the positively oriented 
rectangular path ${\Gamma }_{\varepsilon }$ in $\mathrm{R}$, which consecutively joins the 
vertices $(h, \varepsilon )$, $(h+1, \varepsilon )$, $(h+1, - \varepsilon )$, $(h, - \varepsilon )$, and 
$(h, \varepsilon)$, where 
$\varepsilon $ is chosen sufficiently small and positive so that ${\Gamma }_{\varepsilon }$ lies in $\mathrm{R}$ 
and $-1 < h < 1$. Then
\begin{align}
\int_{{\Gamma }_{\varepsilon }} \dee I & = \big( I(h+1, \varepsilon ) - I(h, \varepsilon ) \big) +
\big( I(h+1, -\varepsilon ) - I(h+1, \varepsilon ) \big) \notag \\
&\hspace{.25in} +\big( I(h, - \varepsilon ) - I(h+1, -\varepsilon ) \big) + \big(I(h, \varepsilon ) - I(h, -\varepsilon ) \big) 
=0. \notag 
\end{align}
This shows that $I$ is single valued on $R$, because ${\Gamma }_{\varepsilon }$ generates the fundamental 
group of $\mathrm{R}$. \medskip  

\noindent That ${\mathcal{A}}_1$ is locally a real analytic function on 
$\mathrm{R}\setminus \{ \ell = 0 \} $ follows from (\ref{eq-s2ss4fsixvnw}) and the fact that $\widetilde{T}$, 
$\widetilde{\vartheta} $, and $I$ are locally real analytic functions
on $\mathrm{R}\setminus \{ \ell = 0 \} $. To show that ${\mathcal{A}}_1$ is continuous on 
$\mathrm{R} \cap \{ \ell = 0 \} $ it suffices to observe that 
$\lim_{\ell \rightarrow 0}\ell \widetilde{\Theta }(h, \ell ) =0$ 
for all $(h, \ell ) \in \mathrm{R} \setminus \{ \ell = 0 \} $. Continuity on $\mathrm{R}$ follows from 
(\ref{eq-s2ss4fsixvnw}), because $\widetilde{T}$ and $I$ are continuous there. From (\ref{eq-s2ss4sixstarnw}) we see 
that $\frac{\partial {\mathcal{A}}_1}{\partial \ell } = -\widetilde{\Theta }$. Thus the assertions about 
$\frac{\partial {\mathcal{A}}_1}{\partial \ell }$ follow from the properties of the function $\widetilde{\Theta }$. \medskip

\noindent 1. The fact that ${\mathcal{A}}_1$ has a real analytic extension to $\partial \overline{\mathrm{R}}$ 
follows from the proof in proposition 2.3. For $(h(s), \pm \ell (s)) \in \partial \overline{\mathrm{R}}$ for 
some $s \in [-1,0)$ we have ${\mathcal{A}}_1(h(s), \pm \ell (s))  = 0$. 
Therefore ${\mathcal{A}}_1$ is nonnegative on $\overline{\mathrm{R}} = 
\mathrm{R} \cup \partial \overline{\mathrm{R}}$. 
Since $(-1,0) \in \partial \overline{\mathrm{R}}$, we get ${\mathcal{A}}_1(-1,0)= 0$. 
The next computation shows that ${\mathcal{A}}_1(1,0) = 4/\pi $. 
\begin{align}
{\mathcal{A}}_1(1,0) & = \frac{1}{\pi } \int^{1}_{-1} \frac{\sqrt{2(1-x)(1-x^2)}}{1-x^2} \, \dee x 
= \frac{\sqrt{2}}{\pi }\, \int^1_{-1} \frac{1}{\sqrt{1+x}} \, \dee x  = \frac{4}{\pi }. \notag 
\end{align}
Thus ${\mathcal{A}}_1$ has a continuous extension to $\overline{\mathrm{R}}$. \medskip 

We now show that $\lim_{h\nearrow \infty}{\mathcal{A}}_1(h,0) = \infty$. First we demonstrate that 
$\lim_{h\nearrow \infty}h\widetilde{T}(h,0)$ $ = \infty$. \medskip 

\noindent \textbf{Proof}. When $h > 1$ and $x \in [-1,1]$ we have $2(h-x)(1-x^2) < 4(h-x)(1-x)$. So 
\begin{displaymath}
\widetilde{T}(h,0) = \frac{1}{\pi } \int^1_{-1} \frac{1}{\sqrt{2(h-x)(1-x^2)}} \dee x > 
\frac{1}{2\pi} \int^1_{-1} \frac{1}{\sqrt{(h-x)(1-x)}} \dee x = I .
\end{displaymath}
Making the successive changes of variables $x = h - u^2$, $v = (h-1)^{-1/2}u$ and $v = \sec s$, we get
\begin{align}
I & = \frac{1}{\pi } \int^{\sqrt{h+1}}_{\sqrt{h-1}} \frac{1}{\sqrt{u^2-(h-1)}} \dee u  
= \frac{1}{\pi } \int^{\sqrt{\frac{h+1}{h-1}}}_1 \frac{1}{\sqrt{v^2-1}} \dee v \notag \\
& = \frac{1}{\pi }\int^{{\sec}^{-1}\sqrt{\frac{1+1/h}{1-1/h}}}_0 \sec s \, \dee s = K. \notag 
\end{align}
Evaluating $K$ we get  
\begin{align}
K & = \frac{1}{\pi } \ln | \sec s + \sqrt{{\sec }^2 s -1}|^{{\sec}^{-1}\sqrt{\frac{1+1/h}{1-1/h}}}_0 \notag \\
& = -\frac{1}{\pi } \ln \sqrt{1-1/h} + \frac{1}{\pi } \ln \big( \sqrt{1+1/h} +\sqrt{2/h} \big) | \notag \\
&= -\frac{1}{2\pi}\ln (1-1/h) +\frac{1}{2\pi } \ln (1+1/h) +\frac{1}{\pi }\ln \big( 1 +\sqrt{\frac{2}{h}(1+1/h)^{-1}}\big) 
\notag \\
& = \frac{1}{2\pi h} + \frac{1}{2\pi h} +\frac{\sqrt{2}}{\pi }h^{-1/2} + \mathrm{O}(h^{-3/2}), \, \, \, 
\mbox{as $h \nearrow \infty$.} \notag
\end{align}
Thus $h \widetilde{T}(h,0) > \frac{1}{\pi } +\frac{\sqrt{2}}{\pi }h^{1/2} + \mathrm{O}(h^{-1/2})$ as $h \nearrow \infty $. 
So $\lim_{h\nearrow \infty} h\widetilde{T}(h,0) = \infty$.  
\vspace{-.05in}
\smallspace
Next we show that $\lim_{h\nearrow \infty}I(h,0) =0$. Now 
\begin{align}
|I(h,0)| & \le \frac{1}{\pi } \int^1_{-1} \frac{|x|}{\sqrt{2(h-x)(1-x^2)}} \dee x \notag \\
& = \frac{1}{\pi }\int^1_0 \frac{x}{\sqrt{2(h-x)(1-x^2)}} \dee x + 
\frac{1}{\pi }\int^0_{-1} \frac{|x|}{\sqrt{2(h-x)(1-x^2)}} \dee x \notag \\
& \le  \frac{1}{\pi }\int^1_0 \frac{x}{\sqrt{2(h-x)(1-x^2)}} \dee x + 
\frac{1}{\pi }\int^1_0 \frac{x}{\sqrt{2(h+x)(1-x^2)}} \dee x \notag \\
& \le \frac{1}{\pi } \frac{1}{\sqrt{h-1}} \int^1_0 \frac{x}{\sqrt{1-x^2}} \dee x + 
\frac{1}{\pi } \frac{1}{\sqrt{h}} \int^1_0 \frac{x}{\sqrt{1-x^2}} \dee x, \notag \\
&\hspace{.5in}\parbox[t]{3.5in}{since $h-x \ge h-1$ and 
$h+x \ge h$ when $x \in [0,1]$} \notag \\
& \le \frac{1}{\pi } \Big( \frac{1}{\sqrt{h-1}} + \frac{1}{\sqrt{h}} \Big), \, \, \mbox{because $ \int^1_0 \frac{x}{\sqrt{1-x^2}} \dee x =1$.} \notag 
\end{align}
So as $h \nearrow \infty $ we see that $|I(h,0)| \searrow 0$. Thus $\lim_{h\nearrow \infty}I(h,0) =0$. \medskip 

Therefore as $h \nearrow \infty$ it follows that ${\mathcal{A}}_1(h,0) = 2h\widetilde{T}(h,0) - I(h,0)$ converges to 
$\infty $. \hfill $\square $ \vspace{-.095in} \medskip  

\noindent 2. Since $\frac{\partial {\mathcal{A}}_1}{\partial h} = \widetilde{T}$ and $\widetilde{T} > 0$ on $\mathrm{R}$, 
it follows that ${\mathcal{A}}_1|_{\overline{\mathrm{R}} \cap \{ \ell = \pm \ell (s_0) \} }$ is strictly increasing. Moreover, 
$\frac{{\partial }^2 {\mathcal{A}}_1}{{\partial h}^2} = \frac{\partial \widetilde{T}}{\partial h} < 0$ on $\mathrm{R}$. To see this we compute
\begin{align}
\frac{\partial \widetilde{T}}{\partial h} & = \frac{1}{\pi } 
\int^{x_{+}}_{x_{-}}  \frac{\partial }{\partial h }  \Big( \big( 2(h-x)(1-x^2) - {\ell }^2 \big)^{-1/2} \Big) \dee x \notag \\
& = -\frac{1}{\pi } \int^{x_{+}}_{x_{-}} (1-x^2)\big( 2(h-x)(1-x^2) - {\ell }^2 \big)^{-3/2} \dee x < 0, \notag 
\end{align}
since the integrand is positive. So the graph of ${\mathcal{A}}_1|_{\overline{\mathrm{R}} \cap \{ \ell = \pm \ell (s_0)  \} }$ 
is strictly convex. Because of convexity, the function  
${\mathcal{A}}_1|_{\overline{\mathrm{R}} \cap \{ \ell = \pm \ell (s_0)  \} }$ is proper. 
Thus if its image were compact, then so 
would be its domain. But this is a contradiction, since its domain $\left[0, \infty ) \right. $  is unbounded. 
Because ${\mathcal{A}}_1(h(s_0), \pm \ell (s_0) ) =0$, the image of ${\mathcal{A}}_1|_{\overline{\mathrm{R}} \cap 
\{ \ell = \pm \ell (s_0) \} }$ is $\left[ 0, \infty ) \right. $. \medskip 

\noindent 3. Let $a >0$. From (\ref{eq-s2ss4sixstarnw}) we find that 
$\frac{\partial {\mathcal{A}}_1}{\partial h}(h,0) = \widetilde{T}(h,0) > 0 $, when $-1 < h <1$ or $h>1$. But 
${\mathcal{A}}_1(-1,0) = 0$, ${\mathcal{A}}_1(1,0) = 4/\pi $, and $\lim_{h \nearrow \infty} {\mathcal{A}}_1(h,0) = \infty$. Thus ${\mathcal{A}}_1([-1,1],0) =[0, 4/\pi ]$ and ${\mathcal{A}}_1([1, \infty ),0) = [4/\pi , \infty )$. Consequently, the 
$a$-level set of ${\mathcal{A}}_1$ intersects the $h$-axis in $\overline{\mathrm{R}}$.  \medskip 

\noindent From (\ref{eq-s2ss4sixstarnw}) and the fact that $\widetilde{T} > 0$ on $\mathrm{R}$, it follows that 
we have $\frac{\partial {\mathcal{A}}_1}{\partial h} > 0$ on ${\mathrm{R}}\cap \{ \ell = 0 \} $. Therefore by the implicit function theorem, near $(h_0, {\ell }_0) \in {\mathcal{A}}^{-1}_1(a)$ there is a real analytic function $\ell \mapsto {\widehat{A}}_a(\ell )$ with $h_0 = {\widehat{A}}_a({\ell }_0)$ such that 
\begin{equation}
a = {\mathcal{A}}_1( {\widehat{A}}_a(\ell ), \ell ) . 
\label{eq-s3ss9onenw}
\end{equation}
Differentiating (\ref{eq-s3ss9onenw}) with respect to $\ell $ and then evaluating the result at $(h_0, {\ell }_0)$ gives 
\begin{equation}
{\widehat{A}}^{\prime }_a( {\ell }_0) = -\frac{\frac{\partial {\mathcal{A}}_1}{\partial \ell }(h_0, {\ell }_0)}
{\frac{\partial {\mathcal{A}}_1}{\partial h }(h_0, {\ell }_0)} = \frac{\widetilde{\Theta }(h_0, {\ell }_0)}
{\widetilde{T}(h_0, {\ell }_0)}, 
\label{eq-s3ss9twonw}
\end{equation}
using (\ref{eq-s2ss4sixstarnw}). ${\widehat{A}}^{\prime }_a$ has a jump discontinuity at 
$\mathrm{R} \cap \{ \ell =0 \}$, because 
\begin{displaymath} 
\lim_{{\ell }_0 \nearrow 0}{\widehat{A}}^{\prime }_a({\ell }_0) 
= \left\{ \begin{array}{rl}
-\frac{1}{2\widetilde{T}(h_0, {\ell }_0)}, & \mbox{if $-1 < h_0 <1$} \\
\smallrowspace -\frac{1}{\widetilde{T}(h_0, {\ell }_0)}, & \mbox{if $h_0 >1$;}
\end{array} \right. 
\end{displaymath}
whereas 
\begin{displaymath} 
\lim_{{\ell }_0 \searrow 0}{\widehat{A}}^{\prime }_a({\ell }_0) 
= \left\{ \begin{array}{rl}
\frac{1}{2\widetilde{T}(h_0, {\ell }_0)}, & \mbox{if $-1 < h_0 <1$} \\
\smallrowspace \frac{1}{\widetilde{T}(h_0, {\ell }_0)}, & \mbox{if $h_0 >1$.}
\end{array} \right. 
\end{displaymath}
Since $\widetilde{T} > 0 $ and $\widetilde{\Theta }(h, \ell ) \in ${\tiny $\left\{ \begin{array}{cl} 
(-1, -\onehalf ) , & \mbox{if $(h, \ell ) \in \mathrm{R} \cap \{ \ell < 0 \} $} \\ 
(\onehalf , 1) , & \mbox{if $(h, \ell ) \in \mathrm{R} \cap \{ \ell > 0 \} $,} \end{array} \right. $} 
it follows that ${\widetilde{A}}^{\prime }_a <0$ on ${\mathrm{R}}\cap \{ \ell < 0 \}$; whereas 
${\widetilde{A}}^{\prime }_a >0$ on ${\mathrm{R}}\cap \{ \ell > 0 \}$. Note that 
${\widehat{A}}^{\prime }_{4/\pi }(0) = 0$. 
Thus a connected component of ${\mathcal{A}}^{-1}_1(a)$ in $\overline{\mathrm{R}}$ is the graph of a 
piecewise real analytic function of $\ell $. Suppose that the domain of ${\widehat{A}}_a$ is the 
compact interval $[{\ell }_{\ast }, {\ell }^{\ast }]$. If ${\ell }^{\ast } <0$, then the function ${\widehat{A}}_a$ would have 
a minimum value $h^{\ast }$ at ${\ell }^{\ast }$, since it is strictly decreasing in $\mathrm{R} \cap \{ \ell < 0 \}$. 
Since $(h^{\ast }, {\ell }^{\ast }) \in \mathrm{R} \cap \{ \ell < 0 \}$, we have ${\widehat{A}}^{\prime }_a({\ell }^{\ast }) < 0$. 
Therefore by the implicit function theorem, there is an $0 > \ell  > {\ell }^{\ast }$ close to ${\ell }^{\ast }$ 
such that ${\widehat{A}}_a$ is defined. But this contradicts the hypothesis that $[{\ell }_{\ast }, {\ell }^{\ast }]$ is 
the domain of ${\widehat{A}}_a$. Therefore ${\ell }^{\ast } > 0$. Then the function ${\widehat{A}}_a$ would have 
a maximum value $h^{\ast }$ at ${\ell }^{\ast }$, since it is strictly increasing in $\mathrm{R} \cap \{ \ell > 0 \}$. 
Since $(h^{\ast }, {\ell }^{\ast }) \in \mathrm{R} \cap \{ \ell >0 \}$, we have ${\widehat{A}}^{\prime }_a({\ell }^{\ast }) >0$. 
Therefore by the implicit function theorem, there is an $\ell  > {\ell }^{\ast }$ close to ${\ell }^{\ast }$ 
such that ${\widehat{A}}_a$ is defined. But this contradicts the hypothesis that $[{\ell }_{\ast }, {\ell }^{\ast }]$ is 
the domain of ${\widehat{A}}_a$. Hence ${\ell }^{\ast }$ does not exist. A similar argument shows that 
${\ell }_{\ast }$ does not exist. Thus the domain of ${\widehat{A}}_a$ is $( -\infty, \infty ) $. This implies 
that the image of ${\widehat{A}}_a$ is $\left[ h_{\ast}(a), \infty ) \right. $, where $h_{\ast }(a) = 
\min_{\ell \in (-\infty, \infty )}{\widehat{A}}_a(\ell ) ={\widehat{A}}_a(0)$ To see that the image of ${\widehat{A}}_a$ is 
unbounded, suppose that $\widetilde{h}= \sup_{\ell \in \left[ -{\ell }(s_0), \infty ) \right. } 
\hspace{-3pt}{\widehat{A}}_a(\ell ) < \infty $. 
Then there is an $\widetilde{\ell } \in {\R }_{\ge 0}$ such that 
$(\widetilde{h}, \widetilde{\ell } ) \in \partial \overline{\mathrm{R}}$. Thus the $a$-level set of 
${\mathcal{A}}_1$ lies in $\overline{\mathrm{R}} \cap \{ h \le \widetilde{h} \} $, 
where $|\ell  | \le  \widetilde{\ell }$. But there is an  
${\ell }^{\dagger} $ in the domain of ${\widehat{A}}_a$ such that  ${\ell }^{\dagger }> \widetilde{\ell }$. Then point 
$({\widehat{A}}_a({\ell }^{\dagger }) , {\ell }^{\dagger }) \in {\mathcal{A}}^{-1}_1(a)$ 
does not lie in $\overline{\mathrm{R}} \cap \{ h \le \widetilde{h} \} $. This is a contradiction. Hence 
$\widetilde{h} = \infty$. 
\medskip  

\noindent Suppose that $a \ne 4/\pi $. Then the line segments 
$\{ (h, -\ell (s)) \in \mathrm{R} \setrule s \in [-1,0) \} $ and 
$\{ (h, \ell (s)) \in \mathrm{R} \setrule s \in [-1,0) \} $, which are parallel to the $h$-axis, each 
intersect the graph of ${\widehat{A}}_a: ( -\infty , \infty ) \longmapsto [h_{\ast }(a), \infty ) $ exactly once, 
since ${\mathcal{A}}_1|_{\mathrm{R} \cap \{ \ell = \pm \ell (s) \} }$ is strictly increasing and has range $(h(s), \infty )$. 
Thus the $a$-level set of ${\mathcal{A}}_1$ with $a>0$ and $\ne 4/\pi $ is connected and is the graph of a 
piecewise real analytic function, whose graph intersects the $h$-axis in 
$\overline{\mathrm{R}}$ exactly once at $(h_{\ast}(a),0)$.  \medskip

\noindent We now look at the $4/\pi $-level set of ${\mathcal{A}}_1$. Since ${\mathcal{A}}_1(1,0) = 4/\pi $, 
the $4/\pi $-level set of ${\mathcal{A}}_1$ is 
nonempty. It is the graph of a piecewise real analytic function ${\widehat{A}}_{4/\pi }$ with 
${\widehat{A}}^{\prime }_{4/\pi }(0) =0$, which is strictly decreasing on $\mathrm{R} \cap \{ \ell < 0 \}$ and 
is strictly increasing on $\mathrm{R} \cap \{ \ell > 0 \}$. The graph of ${\widehat{A}}_{4/\pi }$ intersects 
the $h$-axis at $(1,0)$. The domain of ${\widehat{A}}_{4/\pi }$ 
is $(-\infty, \infty )$ and its range is $[1,\infty)$.   \hfill $\square $ \medskip  

\noindent \textbf{Fact 2.4} The action map of the spherical pendulum is
\begin{equation}
\mathcal{A}: \mathrm{R} \subseteq {\R }^2 \rightarrow  {\R }_{> 0} \times \R \subseteq {\R }^2:  
 \big( h, \ell  \big) \mapsto \big( {\mathcal{A}}_1(h, \ell ), {\mathcal{A}}_2(h, \ell ) \big) 
 \label{eq-twostar}
\end{equation}
is a homeomorphism of $\mathrm{R}$ onto $\big( {\R }_{> 0} \times \R \big) \setminus \{ (1,0) \} $, which 
is a real analytic diffeomorphism on $\mathrm{R} \setminus \{ \ell =0 \} $. This homeomorphism 
extends to a homeomorphism of $\overline{\mathrm{R}} \setminus \{ (1,0) \} $ onto 
$({\R }_{\ge 0} \times \R ) \setminus \{ (1,0) \} $, which is a real analytic diffeomorphism 
on $\overline{\mathrm{R}} \setminus \{ \ell =0 \}$.  \medskip 

\noindent \textbf{Proof.} For every $(h, \ell ) \in \overline{\mathrm{R}}\setminus \{ \ell = 0 \} $ we have 
\begin{displaymath}
D\mathcal{A}(h, \ell ) = 
\begin{pmatrix}
\frac{\partial {\mathcal{A}}_1}{\partial h} (h, \ell ) & \frac{\partial {\mathcal{A}}_2}{\partial h} (h, \ell ) \\
\rowspace \frac{\partial {\mathcal{A}}_1}{\partial \ell } (h, \ell ) & \frac{\partial {\mathcal{A}}_2}{\partial \ell } (h, \ell ) 
\end{pmatrix} = 
\begin{pmatrix} 
\rule{8pt}{0pt}\widetilde{T}(h, \ell ) & 0 \\
-\widetilde{\Theta }(h, \ell ) & 1 
\end{pmatrix} .
\end{displaymath}
The second equality follows from (\ref{eq-s2ss4sixstarnw}). Since $\widetilde{T} >0$ on $\overline{\mathrm{R}}$, 
the action map  $\mathcal{A}$ is a local diffeomorphism on $\mathrm{R} \setminus \{ \ell = 0 \}$. The 
map $\mathcal{A}$ is one to one on $\mathrm{R}$, because 
by point 3 of proposition 2.4 every $a$-level set of ${\mathcal{A}}_1$ on $\overline{\mathrm{R}}$ is the 
graph of a continuous function of $\ell $ on $(-\infty , \infty )$. Thus the $a$-level set of ${\mathcal{A}}_1$ 
intersects the $b\, $-level set of ${\mathcal{A}}_2$ at exactly one point for every $a \ge 0$ and every $b$. 
From the fact that $\mathcal{A}(1,0) = (4/\pi ,0)$ we obtain that 
$\mathcal{A}(\overline{\mathrm{R}} \setminus \{ (1,0) \}) = ({\R }_{\ge 0} \times \R ) \setminus \{ (1,0) \} $. The statement about the extension follows because $\mathcal{A}$ extends to a real analytic mapping on 
$\overline{\mathrm{R}}\setminus \{ 0 \}$ and 
because $\mathcal{A}\big( h(s), \pm \ell (s) \big) = \big( 0, \pm (-s)^{-1/2}(1-s^2) \big) $ 
is a diffeomorphism of $\partial \overline{\mathrm{R}} $ onto $\{ 0 \} \times \R $. \hfill $\square $ \medskip

\subsection{The period lattice and its degeneration}
\label{sec2subsec7} 
%%%%%%%%%%%%

Here we discuss the period lattice of the $2$-torus $T^2_{h, \ell } = {\mathcal{EM}}^{-1}(h, \ell )$ 
when $(h, \ell ) \in \mathrm{R}$ and study its degeneration as $(h, \ell )$ converges to a boundary point 
of $\overline{\mathrm{R}}$. \medskip 

For $(h, \ell ) \in \mathrm{R}$ consider the ${\R}^2$-action on the $2$-torus $T^2_{h, \ell }$ defined by 
\begin{equation}
\Psi : {\R}^2 \times T^2_{h, \ell } \rightarrow T^2_{h, \ell }: (%
\big( (t_1, t_2), p \big) \mapsto {\varphi }^H_{t_1}(p) \comp {\varphi }^L_{t_2}(p).  
\label{eq-s3ss6one}
\end{equation}
Then $\Psi $ is locally transitive at $p \in T^2_{h, \ell }$ because 
$T_pT^2_{h, \ell } = \mathop{\rm span}\nolimits \{ X_H (p), \, X_L (p) \} $.
Since $T^2_{h, \ell }$ is connected, it follows that the action $\Psi $ is
transitive on $T^2_{h, \ell }$. Let ${\mathcal{P}}_{h, \ell } = \{ (t_1, t_2) \in {\R}^2 \, \mathop{\rule[-4pt]{.5pt}{13pt}\, }\nolimits \, {\Psi }_{(t_1,t_2)}(p) = p \} $ be the isotropy group of the ${\R}^2$-action $\Psi $ at $p$. Then $T^2_{h, \ell } = {\R}^2/{\mathcal{P}}_{\ell ,h}$.
Because ${\R}^2$ is abelian, the isotropy group ${\mathcal{P}}_{h, \ell }$ does
not depend on the choice of $p$. Since ${\mathcal{P}}_{h, \ell }$ is closed subgroup of
the Lie group ${\R}^2$, it is a Lie group. Suppose that $\dim {\mathcal{P}}_{h, \ell } \ge 1$. 
Then ${\mathcal{P}}_{h, \ell }$ has a one parameter subgroup, namely, $s
\mapsto (a\, s, b\, s)$ for some $(a,b) \in {\mathbb{R} }^2 \setminus \{(0,0) \}$. 
So $p = {\varphi }^H_{a\, s} \comp {\varphi }^L_{b\, s} (p)$, which implies that 
\begin{equation}
0 = \mbox{${\displaystyle \ttfrac{\dee }{\dee s}}
\rule[-10pt]{.5pt}{20pt} \raisebox{-10pt}{$\, {\scriptstyle s=0}$}$} \hspace{%
-8pt}{\varphi }^H_{a\, s} \comp {\varphi }^L_{b\, s} = a\, X_H(p) + b\, X_L(p).  
\label{eq-s3ss6two}
\end{equation}
But $X_H(p)$ and $X_L(p)$ are linearly independent in $T_pT^2_{h, \ell }$. Hence 
(\ref{eq-s3ss6two}) implies that $a = b =0$, which contradicts the definition
of $a$ and $b$. Thus ${\mathcal{P}}_{h, \ell }$ is a $0$-dimensional Lie group and hence
is discrete. This shows that ${\mathcal{P}}_{h, \ell }$ is a $\Z$-lattice,
called the period lattice. By definition of the functions $T$ and 
$\Theta $, we see that the vectors $\big\{ \big( T(h,\ell ), - \Theta (h, \ell ) \big)^t, \, 
\big( 0 , 2\pi \big)^t \big\} $ form a $\mathbb{Z}$-basis of the period lattice ${\mathcal{P}}_{h, \ell }$. The following calculation serves as a check. Let $(n,m) \in {\Z}^2$.%
Then for $p \in T^2_{h, \ell }$ 
\begin{align}
{\Psi }_{m( T, - \Theta ) + n ( 0 , 2\pi ) }(p) & = 
{\Psi }_{( mT , 2\pi n-m\Theta )}(p) = 
({\varphi }^H_T \raisebox{0pt}{$\scriptstyle\circ \, $} 
{\varphi }^L_{-\Theta })^m \raisebox{0pt}{$\scriptstyle\circ \, $} 
({\varphi }^L_{2\pi })^n(p) = p,  \notag
\end{align}
since ${\varphi }^L_t$ and ${\varphi }^H_{\widetilde{T}\, t} 
\raisebox{0pt}{$\scriptstyle\circ \, $} {\varphi }^L_{-\widetilde{\Theta }\, t}$ are periodic 
of period $2\pi $ on $T^2_{h,\ell }$. Note that the map 
\begin{displaymath}
M: \mathrm{R} \subseteq {\R}^2
\rightarrow \Gl (2, \mathbb{R} ): (h, \ell ) \mapsto 
\raisebox{2pt}{\tiny $\begin{pmatrix} T(h,\ell ) & 0  \\
- \Theta (h,\ell ) & 2\pi \end{pmatrix}$}, 
\end{displaymath}
which to each $(h, \ell ) \in \mathrm{R}$ assigns the
basis of ${\mathcal{P}}_{h, \ell }$, is locally a real analytic
matrix valued function. \medskip

We now look at what happens to ${\mathcal{P}}_{h,\ell }$ as $(h, \ell ) \in \mathrm{R}$ converges to 
a point $(h(s), \ell (s)) \in \partial \overline{\mathrm{R}}$ for some $s \in [-1,0)$. Let $p_{h, \ell } \in 
T^2_{h, \ell }$ and suppose that $p_{h, \ell }$ converges to $p_s \in {\mathcal{EM}}^{-1}((h(s), \ell (s)) = S^1_s$. 
In fact $p_s$ is a relative equilibrium of the Hamiltonian vector field $X_H$ on $(TS^2, \omega )$. 
In other words, the integral curve of $X_H$ starting at $p_s$ is a $1$-parameter subgroup of the 
$S^1$ symmetry generated by the angular momentum Hamiltonian vector field $X_L$. 
So $X_H(p_s)$ and $X_L(p_s)$ are linearly dependent. Thus the isotropy group ${\mathcal{P}}_{h, \ell }$ 
at $p_{h, \ell }$ of the ${\R }^2$-action (\ref{eq-s3ss6one}) degenerates to the isotropy group 
${\mathcal{P}}_s = {\mathcal{P}}_{h(s), \ell (s)}$ of the $\R $-action 
\begin{displaymath}
{\Psi }|_{(\{ 0 \} \times \R ) \times S^1_s}: (\{ 0 \} \times \R ) \times S^1_s \rightarrow S^1_s: (0, t_2) \mapsto 
{\varphi }^H_0 \comp {\varphi }^L_{t_2}(p_s) =  {\varphi }^L_{t_2}(p_s), 
\end{displaymath} 
which is generated by $(0, 2\pi )$. \medskip %

Another way to see that this degeneration is to use the action functions 
$A_1 = {\overline{\pi }}^{\ast } ({\mathcal{A}}_1)$ and $A_2 = {\overline{\pi }}^{\ast }( {\mathcal{A}}_2)$ on 
${\mathcal{EM}}^{-1}({\overline{\mathrm{R}}}^{\ast })$, where ${\overline{\mathrm{R}}}^{\ast } = 
\overline{\mathrm{R}} \setminus \{(1,0) \} $. Here ${\mathcal{A}}_1$ on $\mathrm{R}$ is the first component of 
the action map (\ref{eq-twostar}) and on $\partial \overline{\mathrm{R}}$ is $0$; while $A_2$ on 
${\overline{\mathrm{R}}}^{\ast }$ is $\ell $. Also 
$\overline{\pi }: {\mathcal{EM}}^{-1}({\overline{\mathrm{R}}}^{\ast }) \rightarrow {\overline{\mathrm{R}}}^{\ast } $ is 
the restriction of the energy momentum mapping $\mathcal{EM}$ to the open set 
$TS^2 \setminus {\mathcal{EM}}^{-1}(1,0) = {\mathcal{EM}}^{-1}({\overline{\mathrm{R}}}^{\ast })$. Consider the 
${\R }^2$-action 
\begin{displaymath}
\overline{\Psi }: {\R }^2 \times {\mathcal{EM}}^{-1}({\overline{\mathrm{R}}}^{\ast }) \rightarrow 
{\mathcal{EM}}^{-1}({\overline{\mathrm{R}}}^{\ast }) : \big( (t_1, t_2), q \big) \mapsto 
{\varphi }^{A_1}_{t_1} \comp {\varphi }^{A_2}_{t_2} (q). 
\end{displaymath}
The period lattice ${\mathcal{P}}_q$, which is the isotropy group of 
the action $\overline{\Psi }$ at $q$, is generated by $\{ (2\pi ,0), (0, 2\pi ) \}$, when $q \in  
{\mathcal{EM}}^{-1}(\mathrm{R})$, because the flows of the vector 
fields $X_{A_1}$ and $X_{A_2}$ are periodic of period $2 \pi $. 
When $q \in {\mathcal{EM}}^{-1}(\partial \overline{\mathrm{R}})$, the period lattice ${\mathcal{P}}_q$ is 
generated by $\{ (0, 2\pi ) \} $, because $X_{A_1}(q) = 0$, while $X_{A_2}$ 
has periodic flow of period $2\pi $ on ${\mathcal{EM}}^{-1}(\partial \overline{\mathrm{R}})$. 

\subsection{Monodromy}
\label{sec2subsec8}
%%%%%%%%%%%%%%

In this subsection we show that the spherical pendulum has monodromy by looking at the 
variation of the period lattice along a homotopically nontrivial loop in $\mathrm{R}$. \medskip 

Let $\Gamma :[0,1] \rightarrow \mathrm{R}:t \mapsto \Gamma (t) = \big( h(t), \ell (t) \big) $ 
be a smooth closed curve in $\mathrm{R}$, which encircles the point $(1,0)$ and
generates the fundamental group of $\mathrm{R}$. We transport the period
lattice along $\Gamma $, namely we look at the curve 
\begin{displaymath} 
{\Gamma }^{\ast }M: [0,1] \rightarrow \Gl (2, \R ): t \mapsto M(\Gamma (t))=
\raisebox{2pt}{\tiny $\begin{pmatrix}
T(h(t),\ell (t) ) &  0 \\
- \Theta (h(t) ,\ell (t) ) & 2\pi \end{pmatrix}$. }
\end{displaymath}
We note that 
\begin{equation*}
({\Gamma }^{\ast }M)(1)  
\raisebox{1.5pt}{ {\tiny $ = \begin{pmatrix} 
T(\Gamma (1)) & 0 \\
- \Theta (\Gamma (1)) & 2\pi \end{pmatrix} = $}}  
\raisebox{1pt}{ {\tiny $\begin{pmatrix} 
T(\Gamma (0)) & 0 \\
- \Theta (\Gamma (0)) + 2\pi & 2\pi \end{pmatrix}  = $}}
\raisebox{1pt}{ {\tiny $\begin{pmatrix} 1 & 0 \\ 1 & 1 \end{pmatrix}$}} \, ({\Gamma }^{\ast } M)(0). 
\end{equation*}
In other words, after transporting the period lattice ${\mathcal{P}}_{\Gamma (0)}$ along 
$\Gamma $, its initial basis $\{ \big( T(\Gamma (0)), - \Theta (\Gamma (0)) \big)^t , 
\, \big( 0 , \, 2\pi \big)^t \}$ at $\Gamma (0)$ becomes
the final basis $\{ \big( T(\Gamma (0)), - \Theta (\Gamma (0)) \big)^t 
+ \big( 0, 2\pi \big)^t, \, \big( 0 , \, 2\pi \big)^t \}$ at $\Gamma (1) =
\Gamma (0)$. The matrix of this linear transformation with respect to the
initial basis is $\mathcal{M}=${\tiny $\begin{pmatrix}1 & 0 \\ 1 & 1 \end{pmatrix} $}. 
Thinking of the $2$-torus $T^2_{\Gamma (t)}$ associated to the period
lattice ${\mathcal{P}}_{\Gamma (t)}$, we have a smooth bundle of $2$-tori over $\Gamma $
whose fiber over $\Gamma (t)$ is $T^2_{\Gamma (t)}$. The gluing map of the
fibers over the end points $\Gamma (0)$ and $\Gamma (1)$ is the linear map
of ${\mathbb{R} }^2$ into itself with matrix $\mathcal{M}$. Since 
$\mathcal{M} \in \Gl (2, \Z)$, it maps the lattice $(2\pi {\Z})^2$ into itself. 
Thus $\mathcal{M}$ is a diffeomorphism of the
affine $2$-torus $T^2 = {\R}^2 /(2\pi {\Z})^2$ into itself.
Recalling that action-angle coordinates identify $T^2_{h, \ell }$ with $T^2$,
we see that $\mathcal{M}$ is the \emph{monodromy} map of the $2$-torus
bundle over $\Gamma $. Different choices of the action functions or of the closed 
curve in the homotopy class of $\Gamma $ lead to a $2$-torus bundle, which is isomorphic to the one 
constructed above. The isomorphism class of the new bundle is determined by the 
conjugacy class of the monodromy map $\mathcal{M}$ in $\SSl (2, \Z)$. Consequently, the spherical
pendulum has no global action-angle coordinates. 

\section{Quantum spherical pendulum}
\label{sec3}
%%%%%%%%%%%%%%%%%%%%

We begin with a brief review of the elements of geometric quantization,
which we use here. For details see \cite{sniatycki80}

\subsection{Elements of geometric quantization}
\label{sec3subsec1}
%%%%%%%%%%%%%%%%%%%

\subsubsection{The prequantization line bundle}
%%%%%%%%%%%%%%%%%%%%

The first step in geometric quantization of the symplectic manifold $%
(P,\omega )$ is the construction of a complex line bundle $\lambda :\mathcal{%
L}\rightarrow P$ with connection $\nabla $ such that the curvature of $%
\nabla $ is $-\frac{1}{2 \pi \hbar }\omega $, where $\hbar $ is Planck's
constant divided by $2\pi .$\footnote{We use do not introduce an additional symbol for Planck's constant.} 
For the spherical pendulum $P=T^{\ast }S^{2}$, and $\omega =d\theta ,$ where $%
\theta $ is the Liouville form of the cotangent bundle of the sphere. Hence, 
$\mathcal{L}$ is a trivial bundle, and we can introduce a global
trivializing section $\sigma _{0}:P\rightarrow \mathcal{L}$ such that, for
every vector field $X$ on $P$, the covariant derivative of $\sigma _{0}$ in
direction $X$ is 
\begin{equation}
\nabla _{X}\sigma _{0}=-i\hbar ^{-1}\theta (X)\otimes \sigma _{0}.
\label{trivializing section}
\end{equation}%
Moreover, we introduce a Hermitian form $\langle \cdot ,\cdot \rangle $ on $%
\mathcal{L}$ such that $\langle \sigma _{0},\sigma _{0}\rangle =1.$ The
choice of the trivializing section $\sigma _{0}:P\rightarrow \mathcal{L}$
gives rise to an identification 
\begin{equation}
\mathbb{C}\times P\rightarrow \mathcal{L}:(z,p)\mapsto z\sigma _{0}(p).
\label{identification}
\end{equation}%
Under this identification, every section $\sigma $ of $\mathcal{L}$
corresponds to a complex-valued function $\psi $ on $P$ such that $\sigma
=\psi \sigma _{0}.$

\subsubsection{Bohr-Sommerfeld conditions}
%%%%%%%%%%%%%%%%%%%%%%

Let $\Gamma :[a,b]\rightarrow P:t\mapsto \gamma (t)$ be a curve in $P$. A
section $\sigma $ of $\mathcal{L}$ is covariantly constant along $\Gamma $
if $\nabla _{\dot{\Gamma}(t)}\sigma =0$ for every $t\in \lbrack a,b]$, where 
$\dot{\Gamma}(t)$ denotes the tangent vector of $\Gamma $ at $t.$ If $\sigma
=\psi \sigma _{0}$, then 
\[
\nabla _{\dot{\Gamma}(t)}\sigma =\nabla _{\dot{\Gamma}(t)}(\psi \sigma )=%
\left[ \frac{{\small d}}{{\small dt}}\psi (\Gamma (t))-i\hbar ^{-1}\theta (%
\dot{\Gamma}(t))\psi (\Gamma (t))\right] \sigma _{0}(\Gamma (t)).
\]%
Thus, $\sigma $ is covariantly constant along $\Gamma $ if and only if $\psi
(\Gamma (t))$ satisfies the differential equation%
\[
\frac{{\small d}}{{\small dt}}\psi (\Gamma (t))-i\hbar ^{-1}\theta (\dot{%
\Gamma}(t))\psi (\Gamma (t))=0.
\]%
If $\psi (\Gamma (t))\neq 0$, we divide by $\psi (\Gamma (t))$ and integrate
to get%
\[
\ln \psi (\Gamma (b))-\ln \psi (\Gamma (a))=i\hbar ^{-1}\int_{a}^{b}\theta (%
\dot{\Gamma}(t))dt=i\hbar ^{-1}\int_{\Gamma }\theta ,
\]%
which is equivalent to 
\[
\psi (\Gamma (b))=\psi (\Gamma (a))\exp \left[ i\hbar ^{-1}\int_{\Gamma
}\theta \right] .
\]%
If $\Gamma $ is a closed curve in $P,$ that is if $\Gamma (b)=\Gamma (a)$,
then $\psi (\Gamma (b))=\psi (\Gamma (a))\neq 0$ only if $\exp \left[ i\hbar
^{-1}\int_{\Gamma }\theta \right] =1$. We have proved the following version
of the Bohr-Sommerfeld conditions. \medskip 

\noindent \textbf{Theorem 3.1}\label{Bohr-Sommerfeld} If $\Gamma $ is a closed curve in $P$ and $\sigma $
is a section of $\mathcal{L}$ which is covariantly constant along $\Gamma $,
then the pull-back of $\sigma $ to $\Gamma $ is identically zero unless 
\[
\hbar ^{-1}\int_{\Gamma }\theta =2\pi n
\]%
for some integer $n.$

\subsubsection{Prequantization operators}
%%%%%%%%%%%%%%%%%%%%%%

For each $f\in C^{\infty }(P)$, the Hamiltonian vector field $X_{f}$ of $f$
generates a local 1-parameter group $\exp tX_{f}$ of local diffeomorphisms
of $P$ that preserve the symplectic form $\omega $. There is a unique lift of 
$X_{f}$ to a vector field $\widehat{X}_{f}$ on the prequantization line
bundle $\mathcal{L}$ such that $\exp t\widehat{X}_{f}$ preserves the
connection on $\mathcal{L}$ \cite{sniatycki80}. The prequantization
operator $\mathbf{P}_{f}$ associated to $f$ is given by 
\begin{equation}
\mathbf{P}_{f}\sigma =i\hbar \dbydt \hspace{-8pt} \exp t\widehat{X}_{f}\comp 
\sigma \comp \exp (-tX_{f})  
\label{prequantization 0}
\end{equation}%
for every $\sigma \in S^{\infty }(\mathcal{L}).$ Direct computation gives%
\begin{equation}
\mathbf{P}_{f}\sigma =(-i\hbar \nabla _{X_{f}}+f)\sigma .
\label{prequantization}
\end{equation}

\noindent \textbf{Theorem 3.2} For every $f_{1},f_{2}\in C^{\infty }(P)$, 
\begin{equation}
\lbrack \mathbf{P}_{f_{1}},\mathbf{P}_{f_{2}}]\sigma =i\hbar \, 
\mathbf{P}_{\{f_{1},f_{2}\}}\sigma ,  
\label{commutation relations}
\end{equation}
where $\{f_{1},f_{2}\}=X_{f_{2}}(f_{1})=-X_{f_{1}}(f_{2})=\omega
(X_{f_{2}},X_{f_{1}})$ is the Poisson bracket of $f_{1}$ and $f_{2}$. \medskip 

\noindent \textbf{Proof.}
\begin{eqnarray*}
\lbrack \mathbf{P}_{f_{1}},\mathbf{P}_{f_{2}}]\sigma &=&[(-i\hbar \, 
\nabla _{X_{f_{1}}}+f_{1}),(-i\hbar \nabla _{X_{f_{2}}}+f_{2})]\sigma \\
&=&-\hbar ^{2}(\nabla _{X_{f_{1}}}\nabla _{X_{f_{2}}}-\nabla
_{X_{f_{2}}}\nabla _{X_{f_{1}}})\sigma -i\hbar \, 
(X_{f_{1}}(f_{2})-X_{f_{2}}(f_{1}))\sigma \\
&=&-\hbar ^{2}(\nabla _{\lbrack X_{f_{1}},X_{f_{2}}]}-\ttfrac{i}{%
{\hbar }}\omega (X_{f_{1}},X_{f_{2}}))\sigma -i\hbar \, 
(X_{f_{1}}(f_{2})-X_{f_{2}}(f_{1}))\sigma \\
&=&\hbar ^{2}\nabla _{X_{\{f_{1},f_{2}\}}}\sigma +i\hbar (\omega
(X_{f_{1}},X_{f_{2}})-X_{f_{1}}(f_{2})+X_{f_{2}}(f_{1}))\sigma
\end{eqnarray*}%
because $[X_{f_{1}},X_{f_{2}}]=-X_{\{f_{1},f_{2}\}}$. Moreover,%
\[
\omega
(X_{f_{1}},X_{f_{2}})-X_{f_{1}}(f_{2})+X_{f_{2}}(f_{1})=\{f_{1},f_{2}\}. 
\]%
Therefore, 
\begin{align}
\lbrack \mathbf{P}_{f_{1}},\mathbf{P}_{f_{2}}] &= \hbar ^{2}\nabla
_{X_{\{f_{1},f_{2}\}}}\sigma +i\hbar \, \{f_{1},f_{2}\}\sigma =i\hbar (-i\hbar
\nabla _{X_{\{f_{1},f_{2}\}}}+\{f_{1},f_{2}\})\sigma  \notag \\
&= i\hbar \,\mathbf{P}_{\{f_{1},f_{2}\}}\sigma \textrm{.} \tag*{$\square $}
\end{align}

We refer to the map 
\[
\mathbf{P}:C^{\infty }(P)\times S^{\infty }(L)\rightarrow S^{\infty
}(L):(f,\sigma )\mapsto \mathbf{P}_{f}\sigma 
\]%
as the prequantization map. Equation (\ref{commutation relations}) implies
that the map $f\mapsto \frac{1}{i\hbar }\mathbf{P}_{f}$ is a
representation of the Poisson algebra of $C^{\infty }(P)$ on the space $%
S^{\infty }(\mathcal{L})$, which we call the prequantization representation.

\subsubsection{Polarization}
%%%%%%%%%%%%%%%%%%

In Dirac's formulation of quantum mechanics, the space of quantum states of
the spherical pendulum consists of functions on the spectrum of the complete
set of commuting observables. This idea can be identified with the modern
theory of representations of $\mathbb{C}^{\ast }$ algebras. In the framework
of geometric quantization, it gave rise to the notion of polarization of $(P,\omega )$, 
which is given by a complex involutive Lagrangian distribution $F\ $on
the phase space $P$. For $P=T^{\ast }S^{2}$, the choice of $F$ determines
the representation of the quantum theory of the spherical pendulum.

Let $F$ be a complex involutive Lagrangian distribution $F\ $on $P$. In
other words, $F$ is a complex distribution on $P$ such that $\dim _{\mathbb{C%
}}F=2n,$ where $\dim P=2n.$ Moreover, $\omega (u,v)= 0$ for every pair of vectors $u,v\in F$
attached at the same point in $P$. In the following we
also allow for polarizations with singularities, that is, smooth
distributions that are Lagrangian on an open dense subset of $P.$
Quantization in the $F$-representation leads to the space of quantum states 
\begin{equation}
S_{F}(\mathcal{L})=\{\sigma :P\rightarrow \mathcal{L}\mid \nabla _{u}\sigma
=0\textrm{ for each }u\in F\}.  \label{SFL}
\end{equation}%
In other words, $S_{F}(\mathcal{L})$ is the space of sections of $\mathcal{L}
$ that are covariantly constant along $F$. \medskip 

\noindent \textbf{Definition 3.1}
The space $C_{F}^{\infty }(P)$ of directly quantizable functions in terms of
a polarization $F$ consists of functions $f\in C^{\infty }(P)$ such that,
the flow of the Hamiltonian vector field $X_{f}$ of $f$ preserves the
polarization $F$. \medskip 

For each $f\in C^{\infty }(P)$, the Hamiltonian vector field $X_{f}$ of $f$
lifts to a unique vector field on the prequantization line bundle that
preserves the connection. Hence, if $f\in C_{F}^{\infty }(P)$, it follows
that $\mathbf{P}_{f}$ action leaves $S_{F}(\mathcal{L})$ invariant. \medskip 

\noindent \textbf{Definition 3.2}
Direct quantization in the $F$-representation is given by restricting the
domain of the prequantization map to $C_{F}^{\infty }(P)\times S_{F}(%
\mathcal{L})$ and its codomain to $S_{F}(\mathcal{L})$. In other words, 
\begin{equation}
\mathbf{Q}:C_{F}^{\infty }(P)\times S_{F}(\mathcal{L})\rightarrow S_{F}(%
\mathcal{L}):(f,\sigma ) \longmapsto \mathbf{Q}_{f}\sigma \equiv \mathbf{P}%
_{f}\sigma .  
\label{QF}
\end{equation}

\noindent Quantization in the $F$-representation of functions that are not
in $C_{F}^{\infty }(P)$ requires additional assumptions. \medskip 

\noindent \textbf{Definition 3.3}
A polarization of $(P,\omega )$ is real, if it is a complexification of a
(real) involutive Lagrangian distribution. In other words, 
\begin{equation}
F=D\otimes \mathbb{C}\textrm{,}  \label{DC}
\end{equation}%
where $D$ is an involutive Lagrangian distribution on $P$.

\noindent In the following we assume \medskip 

\noindent \textbf{Condition}\label{condition}
$D$ is locally spanned by Hamiltonian vector fields. \medskip 

\noindent This condition allows for a generalization to polarizations with
singularity.

\subsection{Schr\"{o}dinger quantization}
\label{sec3subsec2}
%%%%%%%%%%%%%%%%%%%%

Schr\"{o}dinger quantization of the spherical pendulum corresponds to the
real polarization tangent to fibres of the cotangent bundle projection $\pi
:T^{\ast }S^{2}\rightarrow S^{2}$. In other words, the Schr\"{o}dinger
polarization is $\ker T\pi \otimes \mathbb{C}$, where 
\[
\ker T\pi =\{u\in T(T^{\ast }S^{2})\mid T\pi (u)=0\}.
\]%
For every $p\in T^{\ast }S^{2}$ and $w\in T_{p}(T^{\ast }S^{2})$, the
evaluation of the Liouville form $\theta $ on $w$ equals the evaluation of $p
$ on $T\pi (w)$, that is, $\theta (w)=p(T\pi (w))$. Hence, $\theta $ vanishes on 
$\ker T\pi $, which implies that the
extension of $\theta $ to the complexification of $\ker T\pi $ vanishes.
Equation (\ref{trivializing section}) implies that the trivializing section $%
\sigma _{0}$ of $\mathcal{L}$ is covariantly constant along $\ker T\pi $.
Thus, every section $\sigma $ of $L$ that is covariantly constant along $%
\ker T\pi $ is of the form $\sigma =\psi \sigma _{0}$, where $\psi $ is a
complex-valued function on $T^{\ast }S^{2}$ that is constant along $\ker
T\pi $. However, functions on $T^{\ast }S^{2}$ that are constant along $\ker
T\pi $ are pull-backs by $\pi :T^{\ast }S^{2}\rightarrow S^{2}$ of functions
on $S^{2}$. Thus, 
\begin{equation}
S_{\ker T\pi }(\mathcal{L})=\{(\pi ^{\ast }\Psi )\sigma _{0}\mid \Psi \in
C^{\infty }(S^{2})\otimes \mathbb{C}\}.  \label{S1}
\end{equation}%
Equation (\ref{S1}) shows that we may identify the space $S_{\ker T\pi }(%
\mathcal{L})$ of sections of $\mathcal{L}$ that are covariantly constant
along $\ker T\pi $ with the space of smooth complex-valued functions on $%
S^{2}.$

The space $C_{\ker T\pi }^{\infty }(T^{\ast }S^{2})$ of functions on $%
T^{\ast }S^{2}$ that are directly quantizable in terms of the Schr\"{o}%
dinger polarization $\ker T\pi $ consists of functions such that their
Hamiltonian vector fields preserve $\ker T\pi $. It contains pull-backs by $%
\pi $ of smooth functions on $S^{2}.$ Moreover, for every vector field $X$ on $S^{2},$ 
its natural extension of $\check{X}$ to $T^{\ast}S^{2}$ is a Hamiltonian vector field with Hamiltonian 
$\theta (\check{X})$ and it preserves $\ker T\pi $. It can be easily verified that 
\begin{equation}
C_{\ker T\pi }^{\infty }(T^{\ast }S^{2})=\{\pi ^{\ast }\check{f}+\theta (%
\check{X})\mid \check{f}\in C^{\infty }(S^{2})\textrm{ and }X \in 
\mathfrak{X}(S^{2})\},  \label{C1}
\end{equation}%
where $\mathfrak{X}(S^{2})$ denotes the space of smooth vector fields on $%
S^{2}$.

It follows from equations (\ref{QF}) and (\ref{prequantization}) that, for
every $\check{f}\in C^{\infty }(S^{2})$, $X \in \mathfrak{X}(S^{2})$
and $\Psi \in C^{\infty }(S^{2})\otimes \mathbb{C}$, 
\begin{equation}
\mathbf{Q}_{\pi ^{\ast }\check{f}}(\pi ^{\ast }\Psi )\sigma _{0}=\pi
^{\ast }(\check{f}\Psi )\sigma _{0},  \label{Qf}
\end{equation}%
and 
\begin{equation}
\mathbf{Q}_{\theta (\check{X})}(\pi ^{\ast }\Psi )\sigma _{0}=-i\hbar \, 
\pi ^{\ast }(X(\Psi ))\sigma _{0}.  \label{QX}
\end{equation}%
Schr\"{o}dinger quantization of functions on $T^{\ast }S^{2}$, which can be
expressed as polynomials in $\pi ^{\ast }\check{f}$ and $\theta (\check{X})$, assigns 
the corresponding polynomial in $\mathbf{Q}_{\pi ^{\ast }%
\check{f}}$ and $\mathbf{Q}_{\theta (\check{X})}$. In this case, the
result depends on the ordering of the factors.

It is usually postulated that the quantum Hamiltonian is 
\[
\mathbf{Q}_{H}(\pi ^{\ast }\Psi )\sigma _{0}=\pi ^{\ast }\left( \left( -%
\ttfrac{{\hbar }^2}{2}\Delta +\sin \vartheta
\right) \Psi \right) \sigma _{0}, 
\]%
where $\Delta $ denotes the Laplace-Beltrami operator on $S^{2}$. There are
several derivations of this result, but each of them requires additional
assumptions, \cite{blattner 1973}, \cite{sniatycki80}, \cite{ward-volkmer}
It is not a direct consequence of prequantization and polarization.

An additional assumption of the Schr\"{o}dinger theory is that the
scalar product of quantum states $(\pi ^{\ast }\Psi _{1})\sigma _{0}$ and $%
(\pi ^{\ast }\Psi _{2})\sigma _{0}$ is given by 
\begin{equation}
(\Psi _{1}\mid \Psi _{2})=\int_{S^{2}}\overline{\Psi }_{1}\Psi _{2} \, \dee \mu ,
\end{equation}%
where $\dee \mu $ is the area form of $S^{2}$. The demand that quantum
observables are given by self-adjoint operators, requires that the
operators in equation (\ref{QX}) should be symmetrized. If the scalar
product is introduced in geometric quantization in terms of half-forms, then
the operators corresponding to equation (\ref{QX}) are symmetric \cite{sniatycki80}.

\subsection{Bohr-Sommerfeld spectrum}
\label{sec3subsec3}
%%%%%%%%%%%%%%%%%%%

Bohr-Sommerfeld quantization corresponds to the polarization $\ker T\mathcal{%
EM}\otimes \mathbb{C}$ tangent to the fibres of the energy-momentum map $%
\mathcal{EM}:T^{\ast }S^{2}\rightarrow \mathbb{R}^{2}$. The range of the
energy momentum map is stratified, with open dense stratum given by the set $%
R$ of regular values. There are two one dimensional strata ${\mathcal{B}}_{+}$ and 
${\mathcal{B}}_{-}$ corresponding to the minimum of energy for a positive or negative value of
the angular momentum, respectively, and two singular points $(-1,0)$ and $(1,0)$.

As before, quantum states are sections of the prequantization line bundle
that are covariantly constant along the polarization. Since fibres of the
energy momentum map are compact, values $(h,\ell l)$ of the energy and the
angular momentum that are in supports of sections of $\mathcal{L}$, which
are covariantly constant along the polarization $\ker T\mathcal{EM}\otimes 
\mathbb{C}$, are restricted by the Bohr-Sommerfeld conditions; see Theorem %
\ref{Bohr-Sommerfeld}. These conditions can be easily described in terms of
the action variables $A_{1},A_{2},$ where $A_{1}$ is a continuous function
of the Hamiltonian $H$ and the angular momentum $L,$ and $A_{2}=L$, see 
proposition 2.4. Recall that $A_{1}$ is a continuous function on $T^{\ast }S^{2}$. Moreover, 
$(A_1)|_{\mathcal{EM}^{-1}(1,0)}=4/\pi $, $(A_1)|_{\mathcal{EM}^{-1}(-1,0)} =0$, and 
$(A_1)|_{\mathcal{EM}^{-1}({\mathcal{B}}_{\pm })} =0$. 

In order to avoid excessive notation, for $(h,\ell)$ in the range of $\mathcal{EM}$, we
write ${\mathcal{A}}_{1}(h,\ell)$ for the value at $(h,\ell)$\ of the
push-forward of $A_{1}$\ by the energy momentum map $\mathcal{EM}$. In other words, $A_1 = 
{\mathcal{EM}}^{\ast }{\mathcal{A}}_1$.  For
each fixed $\ell \ne 0$, the function $h\longmapsto {\mathcal{A}}_{1}(h,\ell )$, which is defined on 
$[h_{\ell},\infty )$, where $h_{\ell}$\ is the minimum of the Hamiltonian $H$\ on
the level set $L=\ell$, is strictly increasing and has range $[0,\infty )$.
Therefore, for every $a \ge 0$, the equation ${\mathcal{A}}_{1}(h,\ell)=a$\ has a unique
solution $h_{\ell}(a)$ for $h\in \lbrack h_{\ell}, \infty )$\ in terms of $a$\ and $\ell$, see point 2 of proposition 2.4. \medskip 

\noindent \textbf{Definition 3.4} The Bohr-Sommerfeld energy-momentum spectrum of the quantum spherical pendulum 
is the set $\mathfrak{S}$ of $(h,\ell)$ in the image of $\mathcal{EM}$ that satisfy the Bohr-Sommerfeld
conditions 
\begin{equation}
\oint A_1 \, \dee \varphi _{1}= 2\pi n\hbar \textrm{ \ \ and \ \ }%
\oint  A_{2} \, \dee \varphi _{2}= 2\pi m\hbar  
\label{Bohr-Sommerfeld 0}
\end{equation}

If $(h,\ell )\in R$, that is, $(h, \ell )$ is a regular value in the image of $\mathcal{EM}$, then  
\begin{equation}
{\mathcal{A}}_{1}(h, \ell )=n\hbar \textrm{ \ \ and \ \ }A_{2}(h,\ell )=\ell=m\hbar 
\label{Bohr-Sommerfeld 2}
\end{equation}
for some integers $n>0$ and $m$. Hence, ${\mathcal{A}}_{1}(h,m\hbar )=n\hbar $, 
and we can express the energy $h$ in terms of $n$ and $m$ and write $h=h(n,m)$.

For $(h,\ell )\in {\mathcal{B}}_{\pm }$, the fibre $\mathcal{EM}^{-1}(h,\ell )$ is a 1-torus 
$\mathbb{T}_{h,\ell}^{1}$. In this case ${\mathcal{A}}_{1}=0$ and the Bohr-Sommerfeld
conditions applied to $\mathbb{T}_{h,\ell }^{1}$ give $\ell = m\hbar $, where $m\in \mathbb{Z}$. 
In this case, we can use equation ${\mathcal{A}}_{1}(h,m\hbar )=0 $ 
to express the energy $h$ in terms of $n=0$ and $m>0$ and write $h=h(0,m)$. 

The fibre of the energy momentum map over the singular point $(-1,0)$ is a
zero-dimensional torus, $\mathcal{EM}^{-1}(-1,0)=(-1,0)$, and the
Bohr-Sommerfeld conditions are satisfied with $n=0$ and $m=0$. Thus, we may
write $h(0,0)=-1$.

It remains to consider the singular point $(1,0)$. Since $A_{1}(1,0)=4/\pi $ and $A_{2}(1,0)=0$, 
equation (\ref{Bohr-Sommerfeld 2}) gives $4/\pi  =n\hbar $, which implies that Planck's constant 
$2\pi \hbar =8/n$. This is a very strong condition on Planck's constant, unlikely to be
satisfied in physics. Therefore, we assume that the Bohr-Sommerfeld
conditions are not satisfied by the unstable equilibrium point $(1,0).$ \medskip 

\noindent \textbf{Conclusion} The Bohr-Sommerfeld energy-momentum spectrum $\mathfrak{S}$ is the range of
a map 
\[
{\mathbb{Z}}_{\ge0} \times \mathbb{Z} \rightarrow \overline{\mathrm{R}} = \mathrm{image}\, \mathcal{EM}:
(n,m) \longmapsto \big( h_m(n),m\hbar \big)
\]%
In other words, 
\begin{equation}
\mathfrak{S=}\left\{ (h_m(n),m\hbar )\in \overline{\mathrm{R}} \setrule 
n\geq 0\right\} . 
 \label{spectrum}
\end{equation}

\noindent \textbf{Definition 3.5}
In physics the pairs $(n,m)\in \mathbb{Z}^{2}$ such that $n\geq 0$ are called 
quantum numbers of the spherical pendulum. \medskip 

For every pair $(n,m)$ of quantum numbers, there exists a non-zero
covariantly constant section $\sigma _{n,m}$ of $\mathcal{L}$ restricted the
fibre $\mathcal{EM}^{-1}(h(n,m),m\hbar ).$ The family of sections 
\begin{equation}
\mathfrak{B}=\left\{ \sigma _{n,m} \setrule (n,m)\in \mathbb{Z}^{2}\textrm{ and }n\geq 0\right\}   
\label{basis}
\end{equation}%
forms a basis in the space of quantum states of the Bohr-Sommerfeld theory.
We consider a Hilbert space $\mathfrak{H}$ of quantum states in which $
\mathfrak{B}$ is an orthonormal basis.
\bigspace
\mbox{}\hspace{.5in}\begin{tabular}{l}
\setlength{\unitlength}{1pt}
\vspace{-.15in} \\
\includegraphics[width=300pt]{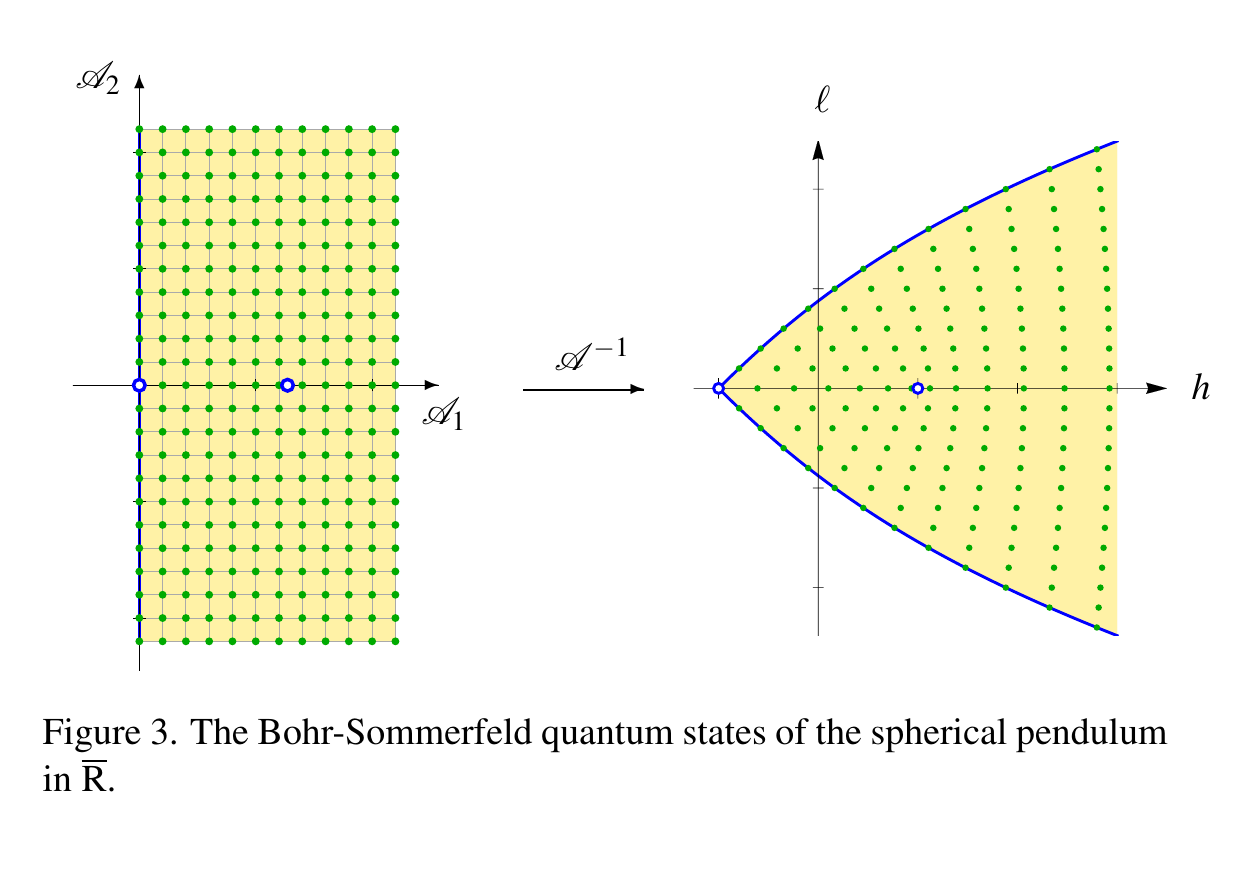} \\
\vspace{-.35in}
\end{tabular} 

According to the general principle of geometric quantization, the space $%
C_{\ker T\mathcal{EM}}^{\infty }(T^{\ast }S^{2})$ of function that are
directly quantizable in the Bohr-Sommer-feld theory consists of smooth
real-valued functions $f$ on $T^{\ast }S^{2}$ such that the Hamiltonian
vector field $X_{f}$ preserves the polarization $\ker T\mathcal{EM}$. \medskip 

\noindent \textbf{Proposition 3.4} If $f\in C^{\infty }(T^{\ast }S^{2})$ is such that $[X_{H},X_{f}]$ and $%
[X_{L},X_{f}]$ are linear combinations of $X_{H}$ and $X_{L},$ then $f$ is a
function of $H$ and $L$. \medskip 

\noindent \textbf{Proof.} Locally, we express $f$ as a function of local angle action coordinates $%
(A_{1},A_{2},\varphi _{1},\varphi _{2})$ and write $f=\breve{f}%
(A_{1},A_{2},\varphi _{1},\varphi _{2}).$ Therefore, 
\begin{align}
X_{f} & = \frac{{\small \partial \breve{f}}}{{\small \partial \varphi }_{1}}X_{\varphi _{1}}+%
\frac{{\small \partial \breve{f}}}{{\small \partial \varphi }_{2}}X_{\varphi _{2}} \notag
\end{align}%
so that 
\begin{align}
\lbrack X_{A_{1}},X_{f}] &= 
\left( X_{A_{1}} \mbox{\small $\frac{\partial \breve{f}}{\partial {\varphi }_{1}}$} \right) 
X_{\varphi _{1}} + \left( X_{A_{1}}  \mbox{\small $\frac{\partial \breve{f}}{\partial {\varphi }_{1}}$ } \right)  
X_{\varphi _{2}}, \notag \\
\lbrack X_{L},X_{f}] &=[X_{A_{2}},X_{f}]=  
\left( X_{A_{2}} \mbox{\small $\frac{\partial \breve{f}}{\partial {\varphi }_{1}}$} \right) 
X_{\varphi _{1}}
+ \left( X_{A_{2}}  \mbox{\small $\frac{\partial \breve{f}}{\partial {\varphi }_{1}}$ } \right)  
X_{\varphi _{2}}. \notag 
\end{align}
Since $X_{A_{1}}$ and $X_{A_{2}}$ are linear combiations of $X_{H}$ and $%
X_{L}$, the assumption that $[X_{H},X_{f}]$ and $[X_{L},X_{f}]$ are linear
combinations of $X_{H}$ and $X_{L}$ implies that 
$X_{A_{i}}\left( \mbox{ \small $\frac{\partial \breve{f}}{\partial {\varphi }_{j}}$} \right) =0$ 
for $i,j=1,2$. Since $X_{A_{i}}=\frac{{\small \partial \breve{f}}}{{\small %
\partial \varphi }_{i}},$ we get $\frac{{\partial }^{2}\breve{f}}{\partial {\varphi }_{i} \partial {\varphi }_{j}}=0$ for $i,j=1,2$.
Integrating the preceding equation gives 
\[
\breve{f}(A_{1},A_{2},\varphi _{1},\varphi _{2})=\check{f}%
_{1}(A_{1},A_{2},\varphi _{1})+\check{f}_{2}(A_{1},A_{2},\varphi _{2}).
\]
Taking into account $\frac{{\small \partial }^{2}{\small \breve{f}}}{{\small %
\partial \varphi }_{1}{\small \partial \varphi }_{1}}=0$ we get that $\check{f}_{1}$ is a linear function of $\varphi _{1}.$ Similarly, $\frac{\small {\partial }^{2}{\small \breve{f}}}{{\small \partial \varphi }_{2}%
{\small \partial \varphi }_{2}}=0$ implies that $\check{f}_{2}$ is a linear
function of $\varphi _{2}$. However, linear functions of angles are not
single-valued. By assumption, $\check{f}$ is single-valued. Therefore, $%
\check{f}$ is independent of $\varphi _{1}$ and $\varphi _{2}$, and
herefore, 
\[
\breve{f}(A_{1},A_{2},\varphi _{1},\varphi _{2})=\check{f}_{3}(A_{1},A_{2}).
\]%
Since $A_{1}$ and $A_{2}$ ar functions of $H$ and $L$, it follows that $f$
restricted to the domain of a local action-angle coordinates $%
(A_{1},A_{2},\varphi _{1},\varphi _{2})$ is a function of the restrictions
of $H$ and $L$ to the same domain.

Domains of local action angle coordinates cover the open set $\mathcal{EM}%
^{-1}(R)\subset T^{\ast }S^{2}$, where $R$ is the set of regular values of
the energy momentum map. Therefore, $f$ restricted to $\mathcal{EM}^{-1}(R)$
is a function of $H$ and $L$ restricted to $\mathcal{EM}^{-1}(R)$. In other
words, $f$ restricted to $\mathcal{EM}^{-1}(R)$ is constant along fibres of
the restriction to $R$ of $\mathcal{EM}:T^{\ast }S^{2}\rightarrow \mathbb{R}%
^{2}$. Since $\mathcal{EM}^{-1}(R)$ is dense in $T^{\ast }S^{2}$ and $f$ is
continuous, it follows that $f$ is constant on fibres of $\mathcal{EM}%
:T^{\ast }S^{2}\rightarrow \overline{\mathrm{R}} \subseteq \mathbb{R}^{2}$. Therefore, there exists a
function $\check{f}:\overline{\mathrm{R}} = \mathrm{image}\, \mathcal{EM}\rightarrow \mathbb{R}$ such
that $f=\check{f}(H,L).$ \hfill $\square $ \medskip 

\noindent \textbf{Corollary 3.5} The space of functions whose Hamiltonian vector fields preserve $\ker T%
\mathcal{EM}$ is%
\[
C_{\ker T\mathcal{EM}}^{\infty }(T^{\ast }S^{2})=\{\check{f}(H,L)\mid \check{%
f}\in C^{\infty }(\mathrm{image}\, \mathcal{EM)}\}\textrm{.} 
\]

Equations (\ref{QF}) and (\ref{prequantization}) imply that the
Bohr-Sommerfeld quantization assigns to a function $\check{f}(H,L)\in
C_{\ker T\mathcal{EM}}^{\infty }(T^{\ast }S^{2})$ an operator $\mathbf{Q}_{\check{f}(H,L)}$ 
on $\mathfrak{H}$ such that operator%
\[
\mathbf{Q}_{\check{f}(H,L)}\sigma _{(h,l)}=\check{f}(h,\ell )\sigma _{(h,\ell )} 
\]%
for every $\sigma _{(h, \ell)}\in \mathfrak{B}$. It should be emphasized that $%
C_{\ker T\mathcal{EM}}^{\infty }(T^{\ast }S^{2})$ does not contain any
function $f$ for which the operator $\mathbf{Q}_{f}$ is not diagonal in
the Bohr-Sommerfeld basis $\mathfrak{B}$.

\subsection{Bohr-Sommerfeld-Heisenberg quantization}
\label{sec3subsec4}
%%%%%%%%%%%%%%%%%%%%%%

\subsubsection{Shifting operators}
%%%%%%%%%%%%%%%%%%

The weakness of the Bohr-Sommerfeld theory is that it does not provide
operators corresponding to transitions between different states. Of course,
there are such operators acting on $\mathfrak{H}$. Since the Bohr-Sommerfeld
basis $\mathfrak{B}$ in $\mathfrak{H}$ has a lattice structure given by
equation (\ref{basis}), there exist shifting operators $\mathbf{a}_{1}$
and $\mathbf{a}_{2}$ on $\mathfrak{H}$ such that%
\begin{equation}
\begin{array}{rl}
\mathbf{a}_{1}\sigma _{n,m} &=\left\{ 
\begin{array}{ccc}
\sigma _{n-1,m} & \textrm{{\small for}} & n\geq 0 \\ 
0 & \textrm{{\small for}} & n=0,%
\end{array} \right.  \\   
\smallrowspace \mathbf{a}_{2}\sigma _{n,m} &=\sigma _{n,m-1}.  \end{array}  
\label{shifting*} 
\end{equation}

\noindent \textbf{Proposition 3.6}
The quantized actions and the shifting operators satisfy the commutation
relations 
\begin{equation}
\begin{array}{rl}
\lbrack \mathbf{Q}_{A_{j}},\mathbf{Q}_{A_{k}}] &=[\mathbf{a}_{j}%
\textrm{,}\mathbf{a}_{k}]=0,  \\
\smallrowspace \lbrack \mathbf{Q}_{A_{j}},\mathbf{a}_{k}] &=-\hbar \, \mathbf{a}_{j}\delta _{jk}, 
\end{array}  
\label{shiftingstar}
\end{equation}
for $j,k=1,2.$

\noindent \textbf{Proof.}
For every basic vector $\sigma _{n,m}$ with $n>0,$%
\begin{eqnarray*}
\lbrack \mathbf{Q}_{A_{1}},\mathbf{a}_{1}]\sigma _{n,m} &=&%
\mathbf{Q}_{A_{1}}\mathbf{a}_{1}\sigma _{n,m}-\mathbf{a}_{1}%
\mathbf{Q}_{A_{1}}\sigma _{n,m}=\mathbf{Q}_{A_{1}}\sigma _{n-1,m}-%
\mathbf{a}_{1}n\hbar \, \sigma _{n,m} \\
&=&(n-1)\hbar \, \sigma _{n-1,m}-n\hbar \, \sigma _{n-1,m}=-\hbar \, \sigma _{n-1,m}
=-\hbar \, \mathbf{a}_{1}\sigma _{n,m}.
\end{eqnarray*}%
Moreover,%
\[
\lbrack \mathbf{Q}_{A_{1}},\mathbf{a}_{1}]\sigma _{0,m}=\mathbf{Q%
}_{A_{1}}\mathbf{a}_{1}\sigma _{0,m}-\mathbf{a}_{1}\mathbf{Q}%
_{A_{1}}\sigma _{0,m}=0=-\hbar \, \mathbf{a}_{1}\sigma _{0,m}.
\]%
Thus, $[\mathbf{Q}_{A_{1}},\mathbf{a}_{1}]=-\hbar \, \mathbf{a}_{1}.
$ Similarly, $[\mathbf{Q}_{A_{2}},\mathbf{a}_{2}]=-\hbar \, \mathbf{%
a}_{2}$. On the other hand, for $n>0$, 
\begin{eqnarray*}
\lbrack \mathbf{Q}_{A_{2}},\mathbf{a}_{1}]\sigma _{n,m} &=&%
\mathbf{Q}_{A_{2}}\mathbf{a}_{1}\sigma _{n,m}-\mathbf{a}_{1}%
\mathbf{Q}_{A_{2}}\sigma _{n,m}=\mathbf{Q}_{A_{2}}\sigma _{n-1,m}-%
\mathbf{a}_{1}m\hbar \, \sigma _{n,m} \\
&=&m\hbar \, \sigma _{n-1,m}-m\hbar \, \sigma _{n-1,m}=0,
\end{eqnarray*}%
and 
\[
\lbrack \mathbf{Q}_{A_{2}},\mathbf{a}_{1}]\sigma _{0,m}=\mathbf{Q%
}_{A_{2}}\mathbf{a}_{1}\sigma _{0,m}-\mathbf{a}_{1}\mathbf{Q}%
_{A_{2}}\sigma _{0,m}=-\mathbf{a}_{1}m\hbar \, \sigma _{0,m}=0.
\]%
Similarly, $[\mathbf{Q}_{A_{1}},\mathbf{a}_{2}]=0$. \hfill $\square $ \medskip 

Since $\mathfrak{B}=\{\sigma _{n,m}\mid (n,m)\in \mathbb{Z}^{2},~n\geq 0%
\mathfrak{\}}\ $is an orthonormal basis in $\mathfrak{H}$, the adjoints 
$\mathbf{a}_{1}^{\dagger }$ and $\mathbf{a}_{2}^{\dagger }$ of the
shifting operators are given by 
\begin{equation}
\begin{array}{rl}
\mathbf{a}_{1}^{\dagger }\sigma _{n,m} &=\sigma _{n+1,m} \\
\smallrowspace \mathbf{a}_{2}^{\dagger }\sigma _{n,m} &=\sigma _{n,m+1},  
\end{array} 
\label{Shifting2}
\end{equation}%
for all $m$ and $n\geq 0$. Taking the adjoints of equations (\ref{shiftingstar}), we get 
\begin{equation}
\begin{array}{rl}
\lbrack \mathbf{a}_{j}^{\dagger }\textrm{,}\mathbf{a}_{k}^{\dagger }]
&=0 \\
\smallrowspace \lbrack \mathbf{Q}_{A_{j}},\mathbf{a}_{k}^{\dagger }] &=\hbar \, \mathbf{a}_{j}^{\dagger }\delta _{jk}. 
\end{array}
\label{Shifting3}
\end{equation}

It follows from the definition, equations (\ref{shifting*}), that the
operators $\mathbf{a}_{1}$ and $\mathbf{a}_{2}$ are analogues of the
lowering operators in the Fock space formulation of field theory and the
states $\sigma _{0,n}$ are ground states for the operator $\mathbf{a}%
_{1} $, \cite{schweber}. Similarly, the operators $\mathbf{a}%
_{1}^{\dagger }$ and $\mathbf{a}_{2}^{\dagger }$ are analogues of the
raising operators.

\subsubsection{Quantization of angles}
%%%%%%%%%%%%%%%%%%%%%%%%%%%

We want to interpret the shifting operators as $\mathbf{a}_{1}$ and $%
\mathbf{a}_{2}$ as quantum operators corresponding to functions $f_{1}$
and $f_{2}$ on the phase space $T^{\ast }S^{2}$ of the spherical pendulum.
In other words, we want to make an identification 
\begin{equation}
\mathbf{a}_{1}=\mathbf{Q}_{f_{1}}\textrm{ \ and \ }\mathbf{a}_{2}=%
\mathbf{Q}_{f_{2}}\textrm{.}  
\label{identification1}
\end{equation}%
Since the shifting operators are not self-adjoint, we cannot expect
functions $f_{1}$ and $f_{2}$ to be real-valued. This means that we have to
extend the Dirac quantization condition 
\begin{equation}
\lbrack \mathbf{Q}_{f_{1}},\mathbf{Q}_{f_{2}}]=i\hbar \, \mathbf{Q}%
_{\{f_{1},f_{2}\}}  
\label{Dirac}
\end{equation}%
to complex-valued functions.

Recall, that the action $A_{1}$ is smooth on the complement $\mathcal{EM}%
^{-1}(R)\setminus L^{-1}(0)$\ of $L^{-1}(0)$ in $\mathcal{EM}^{-1}(R)$,
where $R$ is the set of regular values of the energy momentum map, see proposition 2.3. \medskip 

\noindent \textbf{Proposition 3.6} In $\mathcal{EM}^{-1}(R)\setminus L^{-1}(0)$, the action angle coordinates 
$(A_{1},A_{2}, $ $\varphi _{1},\varphi _{2})$ satisfy satisfy the Poisson bracket
relations%
\begin{equation}
\begin{array}{rl}
\{A_{j},A_{k} \} &= \{ {\mathrm{e}}^{-i\varphi _{j} },{\mathrm{e}}^{-i\varphi _{k} } \}=0 \\
\smallrowspace \{A_{j}, {\mathrm{e}}^{-i\varphi _{k} } \} & =i{\mathrm{e}}^{-i\varphi _{k}}\delta _{jk}. 
\end{array}
\label{commutation 2}
\end{equation}%
for $j,k=1,2$. \medskip 

\noindent \textbf{Proof.}
By definition, the angle action variables satisfy the Poisson commutation
relations%
\[
\{A_{j},A_{k}\}=0\textrm{, \ \ }\{\varphi _{j},\varphi _{k}\}=0\textrm{ \ \ and
\ \ }\{A_{j},\varphi _{k}\}=-\delta _{jk}, 
\]%
where $j,k=1,2$. Hence, 
\[
\{A_{j},{\mathrm{e}}^{-i\varphi _{k}}\}=-i{\mathrm{e}}^{-i\varphi _{k}}\{A_{j},\varphi _{k}\} 
=i {\mathrm{e}}^{-i\varphi _{k}}\delta _{jk}, 
\]%
as required. \hfill $\square $ \medskip 

Comparing equations (\ref{shifting*}) and (\ref{commutation 2}) we see that
that the identification (\ref{identification1}) satisfies the Dirac
quantization condition (\ref{Dirac}) for $f_{1}={\mathrm{e}}^{-i\varphi _{1}}$ and $%
f_{2}={\mathrm{e}}^{-i\varphi _{2}}$ in the open dense subset $\mathcal{EM}%
^{-1}(\mathrm{R})\setminus L^{-1}(0)$ of $T^{\ast }S^{2}$, where $\mathrm{R}$ is the regular
stratum of the range of the energy momentum map $\mathcal{EM}:T^{\ast}S^{2}\rightarrow \mathbb{R}^{2}$. 
In other words, for every $(n,m)\in \mathbb{Z}^{2}$ with $n>0$ and $m\neq 0$, we may set 
\begin{equation}
\mathbf{a}_{1}\sigma _{n,m}=\mathbf{Q}_{{\mathrm{e}}^{-i\varphi _{1}}}\sigma
_{n,m}\textrm{ \ and \ }\mathbf{a}_{2}\sigma _{n,m}=\mathbf{Q}%
_{{\mathrm{e}}^{-i\varphi _{2}}}\sigma _{n,m}.  
\label{identification 2}
\end{equation}

Since ${\mathrm{e}}^{-i\varphi _{1}}$ and ${\mathrm{e}}^{-i\varphi _{2}}$ are well defined on $%
\mathrm{R}\cap L^{-1}(0)$, we can extend the identification (\ref{identification 2}) to $%
\mathrm{R}\cap L^{-1}(0)$ and write 
\begin{equation}
\mathbf{a}_{1}\sigma _{n,0}=\mathbf{Q}_{{\mathrm{e}}^{-i\varphi _{1}}}\sigma
_{n,0}\textrm{ \ and \ }\mathbf{a}_{2}\sigma _{n,0}=\mathbf{Q}%
_{{\mathrm{e}}^{-i\varphi _{2}}}\sigma _{n,0},  
\label{identification 2a}
\end{equation}%
where $n>0$. However, on ${\mathcal{EM}}^{-1}(\mathrm{R}) \setminus L^{-1}(0),$ the Poisson 
brackets involving $A_{1}$ are not defined because $\dee A_{1}$ is not defined there. Hence, the
right hand side of equation (\ref{Dirac}) is not independently defined. This
is a manifestation of the presence of monodromy in the spherical pendulum,
which will be discussed in Section 3.4.4.

\subsubsection{Boundary conditions}
%%%%%%%%%%%%%%%%%%%%%%%

It remains to extend operators $\mathbf{Q}_{{\mathrm{e}}^{-i\varphi _{1}}}$ and 
$\mathbf{Q}_{{\mathrm{e}}^{-i\varphi _{2}}}$ to quantum states $\sigma _{0,m}$
supported on the boundary of $\mathcal{EM}^{-1}(\mathrm{R})$. This is analogous to
extending Schr\"{o}dinger quantization to the cotangent bundle of a manifold
with boundary and corners.

The angle functions are not globally defined. In particular, the functions ${\mathrm{e}}^{-i\varphi _{1}}$ 
and ${\mathrm{e}}^{-i\varphi _{2}}$ are not defined at the singular
points $(0,-1)$ and $(0,1)$. Moreover, ${\mathrm{e}}^{-i\varphi _{1}}$ is not defined
when $A_{1}=0$.  

In order to extend ${\mathrm{e}}^{-i\varphi _{1}}$ to a globally defined function, we
choose a smooth function $\chi _{1}(h,l)$ on $\overline{\mathrm{R}} = \mathrm{image}\, \mathcal{EM}, $
which is identically $1$ on a neighbourhood of $\mathrm{R} \cap \mathfrak{S}$ in $\mathrm{R}$,
where $\mathfrak{S}$ is the Bohr-Sommerfeld energy spectrum (\ref{spectrum}), and vanishes 
to infinite order on the boundary $\partial \mathrm{R}$ of $\mathrm{R}$. The
product $f_{1}=\chi _{1}(L,H){\mathrm{e}}^{-i\varphi _{1}}$ is a globally defined
function on $T^{\ast }S^{2}$ that vanishes to infinite order on $\mathcal{EM}^{-1}(\partial \mathrm{R})$ 
and satisfies the Poisson bracket relations 
\[
\{A_{j},f_{1}\}=\{A_{j},\chi _{1}(L,H){\mathrm{e}} ^{-i\varphi _{1}}\} 
=-i\chi _{1}(L,H){\mathrm{e}}^{-i\varphi _{1}}\{A_{j},\varphi _{1}\} 
=f_{1}\delta _{j1}
\]%
for $j=1,2.$ Since $f_{1}=\chi _{1}(L,H){\mathrm{e}}^{-i\varphi _{1}}$ vanishes to
infinite order on $\partial \mathrm{R} \cap \mathfrak{S}$ and derivatives of $\chi _{1}$ vanish to 
infinite order on $\mathfrak{S}$, it follows that we can make the identification 
\begin{equation}
\mathbf{Q}_{f_{1}}\sigma _{n,m}=\mathbf{Q}_{\chi _{1}{\mathrm{e}}^{-i\varphi _{1}}}\sigma _{n,m}
=\mathbf{a}_{1}\sigma _{n,m}  
\label{identification 3}
\end{equation}%
for all $(n,m)\in \mathbb{Z}^{2}$ with $n>0$ and $m\neq 0$. On the other
hand, $f_{1}=\chi _{1}(L,H){\mathrm{e}}^{-i\varphi _{1}}$ vanishes to infinite order on 
$\partial \mathrm{R} \cap \mathfrak{S}$. Hence, we may set 
\begin{equation}
\mathbf{Q}_{f_{1}}\sigma _{0,m}=0=\mathbf{a}_{1}\sigma _{0,m}
\label{identification 3a}
\end{equation}%
This identification is independent of the choice of $\chi _{1}$ satisfying
the required conditions. In order to keep the notation simple, in the
following we omit $\chi _{1}$ and write 
\begin{equation}
\mathbf{Q}_{{\mathrm{e}}^{-i\varphi _{1}}}\sigma _{n,m} 
=\mathbf{a}_{1}\sigma _{n,m}  
\label{identification 4}
\end{equation}%
for all $(n,m)\in \mathbb{Z}^{2}$ with $n\geq 0$ and $m\neq 0$. It remains
to consider the action of 

The function ${\mathrm{e}}^{-i\varphi _{2}}$ is defined on $\mathcal{EM}^{-1}(\partial
\mathrm{R} \setminus \{(-1,0),(1,0)\}$. As before, we can extend ${\mathrm{e}}^{-i\varphi _{2}}$
to the whole of $T^{\ast }S^{2}$ by multiplying ${\mathrm{e}}^{-i\varphi _{2}}$ by an
appropriate function $\chi _{2}$ of $H$ and $L$. By assumption $(1,0)\notin 
\mathfrak{S}$, so vanishing of $\chi _{2}$ at $(1,0)$ does not affect the
identification of $\mathbf{a}_{2}$ with $\mathbf{Q}_{{\mathrm{e}}^{-i\varphi _{2}}}$. 
However, the point $(-1,0)\in \mathfrak{S}$, corresponding to $(n,m)=(0,0)$, is a corner of 
the range of $\mathcal{EM}$ and $\varphi _{2}$
is not defined on $\mathcal{EM}^{-1}(-1,0)$. Therefore, the identificalion 
\begin{equation}
\mathbf{Q}_{{\mathrm{e}}^{-i\varphi _{2}}}\sigma _{0,0}= 
\mathbf{a}_{2}\sigma _{0,-1}  
\label{identification 5}
\end{equation}%
is an essential extension of the definition of $\mathbf{Q}_{{\mathrm{e}}^{-i\varphi _{2}}}$, 
which does not follow from the Dirac quantization conditions. Making this identification, we write 
\begin{equation}
\mathbf{Q}_{{\mathrm{e}}^{-i\varphi _{2}}}\sigma _{n,m}= 
\mathbf{a}_{2}\sigma _{n,m}  
\label{identification 6}
\end{equation}%
for all $(n,m)\in \mathbb{Z}^{2}$ with $n\geq 0$.

Having identified shifting operators $\mathbf{a}_{1}$ and $\mathbf{a}%
_{2}$ with quantizations of ${\mathrm{e}} ^{-i\varphi _{1}}$ and ${\mathrm{e}}^{-i\varphi _{2}}$,
respectively, we observe that the adjoint operators $\mathbf{a} _{k}^{\dagger }$ may be identified 
with quantization of ${\mathrm{e}}^{i\varphi _{k}}$. In other words,%
\begin{equation}
\mathbf{a}_{k}^{\dagger }=\mathbf{Q}_{{\mathrm{e}}^{i\varphi _{k}}}
\label{identification 7}
\end{equation}%
for $k=1,2.$ Moreover, as in the Schr\"{o}dinger theory, to a function $f$
on $T^{\ast }S^{2}$, which can be expressed as a polynomial in $A_{1},$ $%
A_{2}$, $e^{-i\varphi _{1}}$ and ${\mathrm{e}}^{-i\varphi _{2}}$, we can assign the
corresponding polynomial in $\mathbf{Q}_{A_{1}}$, $\mathbf{Q}_{A_{2}}$, 
$\mathbf{Q}_{{\mathrm{e}}^{-i\varphi _{1}}}$, and 
$\mathbf{Q}_{{\mathrm{e}}^{-i\varphi _{2}}}.$ In this case, the result depends on the ordering of the factors.

\subsection{Quantum monodromy}
%%%%%%%%%%%%%%%%%%%%%%

In this subsection we will discuss \medskip 

1. the definition of quantum monodromy; 
\smallspace
\indent 2. \parbox[t]{4in}{show that the quantum spherical 
pendulum has quantum monodromy;} 
\smallspace
\indent 3. \parbox[t]{4in}{read off the classical monodromy of the spherical 
pendulum from the joint spectrum of its Bohr-Sommerfeld-Heisenberg quantization.} \bigskip 

We begin by defining quantum monodromy. Our discussion leans heavily on the treatment of 
Vu Ngoc \cite{vungoc}. Let $B$ be an open $2$-disk in ${\R }^2$, which is contained in $\mathrm{R}$ 
and is centered at $c \in \mathrm{R}$. The intersection $B \cap \mathfrak{S}$ of the Bohr-Sommerfeld 
spectrum $\mathfrak{S}$ with $B$ is a \emph{local lattice}, because the image of 
$B \cap \mathfrak{S}$ under the homeomorphism given by the action map ${\mathcal{A}}$ (\ref{eq-twostar}) 
is the intersection of the open subset ${\mathcal{A}}(B)$ with the standard lattice 
$2\pi \hbar \, {\Z }^2$. We call the pair $(B, \mathcal{A}|_B)$ a \emph{local chart} at $c$ for the 
Bohr-Sommerfeld spectrum $\mathfrak{S}$. Let $\bigcup_{\alpha \in \mathcal{I}} B_{\alpha }$ be an 
open covering of $\mathrm{R}$ by $2$-disks each centered at $c_{\alpha}$. Suppose that 
$(B_{\alpha }, \mathcal{A}|_{B_{\alpha }})$ and $(B_{\beta }, \mathcal{A}|_{B_{\beta }})$ are local 
charts for $\mathfrak{S}$ and that $c_{\alpha} \in B_{\alpha } \cap B_{\beta }$. From the construction of 
action angle coordinates in section 2.4, it follows that the chart transition map  
\begin{equation}
A_{\alpha \beta }: B_{\alpha } \cap B_{\beta } \subseteq \mathrm{R} \rightarrow \SSl (2, \Z) : 
(h, \ell ) \longmapsto {\mathcal{A}}|_{B_{\beta }} (h, \ell ))\, \comp 
\big( {\mathcal{A}}|_{B_{\alpha }}(h, \ell ) \big)^{-1}
\label{eq-sec3.4.4onenew}
\end{equation}
is locally constant. Let $\tau : \mathcal{L}\rightarrow \mathrm{R}$ be the ${\Z }^2$-bundle over $\mathrm{R}$ 
with local trivialization given by the top row of the commutative diagram \medskip
\begin{displaymath}
\begin{CD}
{\tau }^{-1}_{\alpha }(U_{\alpha }) \subseteq \mathcal{L} @>{\tau }_{\alpha }>> B_{\alpha } \times {\Z}^2 
\subseteq \mathrm{R} \times {\Z }^2 \\
@VV\tau |_{{\tau }^{-1}_{\alpha }(U_{\alpha })}V                  @VV{\pi}_{\alpha }V\\
B_{\alpha} \subseteq {\R}^2 @>\mathrm{id}>> B_{\alpha} \subseteq \mathrm{R} 
\end{CD}
\end{displaymath}
\medskip 

\noindent where the left vertical arrow is the restriction of the bundle projection map to ${\tau }^{-1}_{\alpha }(U_{\alpha })$ and the left vertical arrow is the projection map on the first factor. The bottom horizontal map is 
the identity mapping. More explicitly, the image of the fiber ${\mathcal{L}}_{\widetilde{c}}$ over the 
point $\widetilde{c} \in B_{\alpha }$ under the trivialization mapping ${\tau }_{\alpha }$ is 
$\big( \widetilde{c}, (v^1_{c_{\alpha }}, v^2_{c_{\alpha }}) \big) \in B_{\alpha } \times {\Z }^2$, where 
$\{ v^1_{c_{\alpha }}, \, v^2_{c_{\alpha }} \} $ is a $\Z$-basis of the ${\Z }^2$-lattice ${\mathcal{L}}_{c_{\alpha }}$. 
The local transition maps for the bundle $\mathcal{L}$ are given by 
\begin{equation}
\begin{array}{l}
{\tau }_{\alpha \beta }: (B_{\alpha } \cap B_{\beta }) \times {\Z }^2 \rightarrow 
(B_{\alpha } \cap B_{\beta }) \times {\Z }^2: \\
\smallrowspace \hspace{.35in} (\widetilde{c}, v) \longmapsto \big( \mathcal{A}(\widetilde{c}), 
A_{\alpha \beta }(\widetilde{c})v \big) .
\end{array}
\label{eq-sec3.4.4twonew}
\end{equation}
The isomorphism class\footnote{Note that 
$M = \{ A_{\alpha \beta } \} $ is a $1$ \v{C}ech cocycle with values in $\SSl (2, \Z )$ for the covering 
$\bigcup_{\alpha \in \mathcal{I}} B_{\alpha }$ of $\mathrm{R}$, since 
$A_{\alpha \beta } A_{\beta \alpha } = \mathrm{id}$ and  
$A_{\alpha \beta } A_{\beta \gamma } A_{\gamma \alpha } = \mathrm{id}$, 
for every $\alpha $, $\beta $, 
and $\gamma \in \mathcal{I}$. The \v{C}ech cohomology class $[M]$, associated to the \v{C}ech cocycle $M$, corresponds to the isomorphism class of the bundle $\mathcal{L}$. See \cite[p.40--41]{hirzebruch}.} of the bundle 
$\mathcal{L}$ is the \emph{quantum monodromy} of the Bohr-Sommerfeld spectrum $\mathfrak{S} \subseteq \mathrm{R}$.  \medskip 

Action-angle coordinates provide a local trivialization of the bundle $\pi : {\mathcal{EM}}^{-1}(\mathrm{R}) 
\rightarrow \mathrm{R}$. From their construction it follows that the bundle $\mathcal{L} \rightarrow 
\mathrm{R}$ is isomorphic to the bundle $\mathcal{P} \rightarrow \mathrm{R}$ of period lattices. 
But the bundle $\mathcal{P}$ is not trivial, since the spherical pendulum has classical monodromy. 
Thus the ${\Z }^2$-lattice bundle $\mathcal{L}$ is not trivial. In other words, the quantum 
monodromy of the spherical pendulum is nontrivial. \medskip 

We now give a geometric procedure for determining the quantum monodromy of the 
Bohr-Sommerfeld spectrum $\mathfrak{S}$. Consider the neighboring quantum numbers 
\begin{equation}
(n,m), \, \, \, (n+1, m), \, \, \, (n+1, m+1), \, \, \, \mathrm{and} \, \, \, (n, m+1).
\label{eq-sec3ss4.4zero}
\end{equation}
Their correspoding quantum states ${\sigma }_{n,m}$, ${\sigma }_{n+1,m}$, ${\sigma }_{n+1,m+1}$, and 
${\sigma }_{n,m+1}$ are obtained by applying the shifting operators ${\mathbf{a}}^{\dagger}_1$, 
${\mathbf{a}}^{\dagger}_2{\mathbf{a}}^{\dagger}_1$, and 
${\mathbf{a}}_1{\mathbf{a}}^{\dagger}_2{\mathbf{a}}^{\dagger}_1$, respectively, to ${\sigma }_{n,m}$. The corresponding spectral values in $\mathrm{R}$ of the quantized spherical pendulum, given by 
\begin{equation}
\begin{tabular}{ccc}
$\big( h_{\hbar \, m}( \hbar \, n), \hbar \, m \big) $, & \quad &  
$\big( h_{\hbar \, m }(\hbar (n+1)), \hbar \, m \big) $, \\ 
\rowspace $\big( h_{\hbar (m+1)}(\hbar (n+1)), \hbar (m+1) \big) $ , & \quad &  
$\big( h_{\hbar (m+1)}(\hbar \, n), \hbar (m+1) \big) $,
\end{tabular} 
\label{eq-sec3ss4.4one}%
\end{equation}
form the vertices of a \emph{spectral quadrilateral} $Q_{n,m}$ with $(n,m)$ the lower left hand vertex of 
$Q_{n,m}$. \medskip 

Let $\Gamma $ be a positively oriented, closed, non-self-intersecting polygonal path in $\mathrm{R}$ which \medskip 

1. encircles the point $(1,0)$; 
\smallspace
\indent 2. \parbox[t]{4in}{ passes consecutively throught the vertices 
\begin{displaymath}
(n_1, m_1), \, \, \, (n_2, m_2),  \ldots , (n_{\ell }, m_{\ell }), \, \, \, (n_{\ell +1}, m_{\ell +1}) = (n_1, m_1)
\end{displaymath}
of the spectral quadrilaterals $Q_{n_1,m_1}$, $Q_{n_2,m_2}$, $\ldots Q_{n_{\ell }, m_{\ell }}$, 
$Q_{n_1,m_1}$.} \bigskip 

\noindent \textbf{Claim 3.7} For each $1 \le i \le \ell $ there is a finite shortest sequence ${\mathbf{a}}^i$ 
of shifting operators, each member of which is taken from $\{ {\mathbf{a}}_1, \, {\mathbf{a}}_2, \, 
{\mathbf{a}}^{\dagger}_1, \, {\mathbf{a}}^{\dagger }_2 \} $, that shifts the vertex $(n_i, m_i)$ of 
the spectral quadrilateral $Q_{n_i,m_i}$ to the vertex $(n_{i+1}, m_{i+1})$ of $Q_{n_1+1, m_i+1}$. \medskip 

\noindent \textbf{Corollary 3.8} For each $1 \le i \le \ell $ the image of $Q_{n_i,m_i}$ under the 
operator ${\mathbf{a}}^i$ is $Q_{n_1+1, m_i+1}$. \medskip 

\noindent The \emph{lower left hand corner} of the spectral quadrilateral $Q_{n,m}$ is given by the 
spectral values corresponding to the edge joining $(n,m)$ to $(n+1,m)$ and the spectral values 
correspongind to the edge joining $(n,m)$ to $(n,m+1)$. Because the quantum spectral values are 
determined by the intersection of level sets of the action function ${\mathcal{A}}_1$ with a level set 
of the action function ${\mathcal{A}}_2$, to order 
$\mathrm{O}(h)$ the lower left hand corner of the spectral quadrilateral is given by the row vectors 
of the derivative of the action map $\mathcal{A}$ (\ref{eq-twostar}) at $p_{n,m}= 
(h_{\hbar \, m}(\hbar \, n), \hbar \, m)$. So  
\begin{equation}
D\mathcal{A} (p_{n,m}) = 
\begin{pmatrix}
\frac{\partial {\mathcal{A}}_1}{\partial h}(p_{n,m}) & \frac{\partial {\mathcal{A}}_2}{\partial h}(p_{n,m}) \\
\smallrowspace 
\frac{\partial {\mathcal{A}}_1}{\partial \ell }(p_{n,m}) & \frac{\partial {\mathcal{A}}_2}{\partial \ell }(p_{n,m}) 
\end{pmatrix} = \begin{pmatrix}
\widetilde{T}(p_{n,m}) & 0 \\
\smallrowspace -\widetilde{\Theta }(p_{n,m}) & 1 
\end{pmatrix} , 
\label{eq-sec3ss4.4four}
\end{equation}
using (\ref{eq-s2ss4sixstarnw}) and the fact that ${\mathcal{A}}_2(h, \ell ) = \ell $. \medskip 

We now look at the \emph{variation} of $D\mathcal{A}$ along the lower left hand vertices occuring 
on the polygonal path $\Gamma $ in $\mathrm{R}$. On the one hand because $\widetilde{T}$ is a continuous function on $\mathrm{R}$, its variation along $\Gamma $ is $0$. On the other hand, from fact 2.2 
the variation of $-\widetilde{\Theta }$ along $\Gamma $ is $1$. Thus the variation of $D\mathcal{A}$ 
along $\Gamma $ is {\tiny $\begin{pmatrix} 0 & 0 \\ 1 & 0 \end{pmatrix} $}. Since the column vectors 
of $D\mathcal{A}(p_{n,m})$ form a basis of the period lattice of the $2$-torus 
${\mathbb{T}}^2_{p_{n,m}}$, the monodromy matrix of the classical 
spherical pendulum along $\Gamma $ is the sum of the identity matrix 
{\tiny $\begin{pmatrix} 1 & 0 \\ 0 & 1 \end{pmatrix}$} and the variation, namely, 
{\tiny $\begin{pmatrix} 1 & 0 \\ 1 & 1 \end{pmatrix} $}. See figure 4.
\bigspace
\par \noindent \hspace{.75in} \begin{tabular}{l}
\setlength{\unitlength}{1pt}
\vspace{-.25in}\\
\includegraphics[width=225pt]{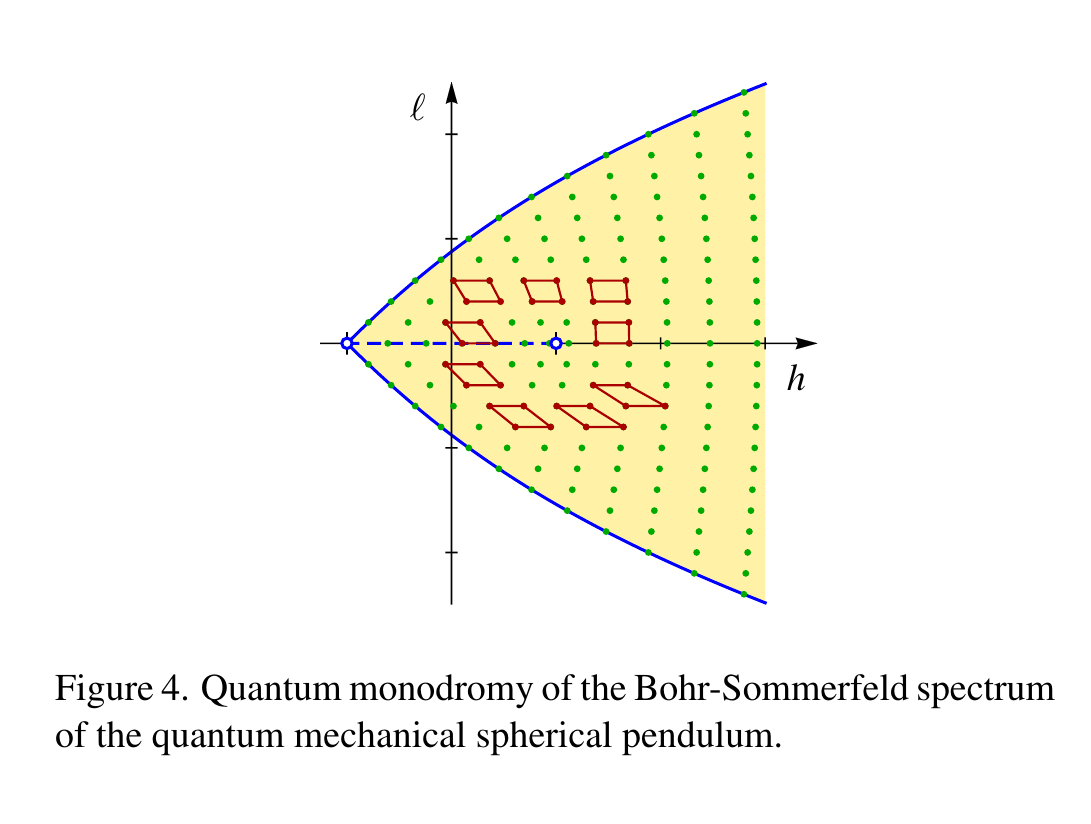}
\end{tabular}

\section{Concluding remarks}

\begin{itemize}
\item We have extended the Bohr-Sommerfeld quantization of the spherical
pendulum to a full quantum theory, which we call the
Bohr-Sommerfeld-Heisenberg (BSH) quantization. Our approach leads to a
matrix formulation of Born, Jordan and Heisenberg \cite{born-jordan}, \cite%
{born-jordan-heisenberg} with quantum operators expressed as matrices in the
Bohr-Sommerfeld basis. According to Mehra and Rechenberg \cite[p.265]{mehra-rechenberg} the connection between quantum mechanics and
the Bohr-Sommerfeld theory of multiply periodic systems \textquotedblleft
seems to have been lost completely in the matrix approach\textquotedblright\
until it was re-established by Wentzel \cite{wentzel}.

\item We had an advantage of being able to rely on the guiding principle of
geometric quantization for an open dense subset of the phase space in which
all our constructions were regular. In treating nowhere dense sets of
singular points, we followed Dirac's Principles of Quantum Mechanics \cite%
{dirac 30} 

\item In our presentation, we have included the geometric quantization
setting of the Schr\"{o}dinger theory. It shows that the main difference
between the BSH theory and the Schrodinger theory is the choice of
polarization. For a completely integrable system, the energy momentum map is
regular in an open dense subset of the phase space. Hence in the BSH
quantization, we have to deal with polarization with singularities on the
boundary of that set. This leads to difficulties analogous to those that
appear in formulating the Schrodinger theory on a singular space.

\item In the case of the spherical pendulum, the Bohr-Sommerfeld energy
spectrum differs from the Schr\"{o}dinger energy spectrum. In the
quasiclassical limit of $\hslash $ close to zero, the Bohr-Sommerfeld
spectrum and the Schr\"{o}dinger spectrum of $\mathbf{Q}_{A_{1}}$ differ
by $\frac{1}{2}\hslash $. 
\item In physics, Planck's constant $2\pi \hslash $ is approximately $%
6.626\times 10^{-34}$ joules. However, in the process of
quantization, $\hslash $ is treated as a parameter. In the Schr\"{o}dinger
quantization, the reperesentation space is independent of $\hslash $, while
the quantum operators depend on $\hslash $ explicitly. In the BSH
quantization, presented here, the representation space $\mathfrak{H}$
depends on $\hslash $ because it is defined in terms of a basis $\mathfrak{B}
$ consisting eigensections supported on fibres of the energy-momentum map
that satisfy the Bohr-Sommerfeld conditions. For every fibre of the
energy-momentum map, there exists a value of $\hslash $, treated as a
parameter, for which this fibre satisfies Bohr-Sommerfeld conditions.
\item For every $\hslash \neq 4/\pi n$, where $n\in \mathbb{N}$, our
construction gives well defined shifting operators ${\mathbf{a}}_{1},$ $%
\mathbf{a}_{2}$ and their adjoints $\mathbf{a}_{1}^{\dagger },$ $%
\mathbf{a}_{2}^{\dagger }$ on the representation space $\mathfrak{H}$.
It is a consequence monodromy that the interpretation of the shifting
operators $\mathbf{a}_{1},$ $\mathbf{a}_{2}$ as quantization of 
${\mathrm{e}}^{-i\varphi _{1}}$ and ${\mathrm{e}}^{-i\varphi _{2}}$ fails to be global even on the
set of regular values of the energy momentum map.

\end{itemize}

%The reference entries Schweber and Wentzel needs to be filled in. 


\begin{thebibliography}{99}
\bibitem{blattner 1973} R.J. Blattner, \textquotedblleft Quantization in
representation theory\textquotedblright\ in \textit{Harmonic Analysis on
Homogeneous Spaces}, edited by E.T. Taam. \textit{Proc. Sym. Pure Math.} vol. 
\textbf{26} pp. 146--165. A.M.S., Providence, R.I. 1973.

\bibitem{blattner 1975} R.J. Blattner, \textquotedblleft Pairing of
half-form spaces\textquotedblright , in \textit{G\'{e}om\'{e}trie
Symplectique et Physique Math\'{e}matique\textquotedblright }, Coll. Int.
C.N.R.S. No 237, pp. 175--186, C.N.R.S., Paris, 1975.

\bibitem{blattner 1977} R.J. Blattner, \textquotedblleft The metalinear geometry of
non-real polarizations\textquotedblright ,\textit{\ Lecture Notes in
Mathematics,} vol. \textbf{570}, pp. 11--45, Springer, New York 1977.

\bibitem{born-jordan} M. Born and P. Jordan, \textquotedblleft Zur
Quantenmechanik\textquotedblright , \textit{Z. f\"{u}r Phys.}, \textbf{34} (1925)
858--888.

\bibitem{born-jordan-heisenberg} M. Born, W Heisenberg and P. Jordan,
\textquotedblleft Zur Quantenmechanik II\textquotedblright , \textit{Z. Phys}., \textbf{35} (1925) 557--615. 

\bibitem{cushman-bates} R.H. Cushman and L.M. Bates, \textit{Global aspects of
classical integrable systems, second edition}, Birkh\"{a}user, Basel, 2015.

\bibitem{cushman-duistermaat} R. Cushman and J.J. Duistermaat,
\textquotedblleft The quantum mechanical spherical
pendulum\textquotedblright , Bull. AMS, 19 (1988) 475--479.

\bibitem{cushman-sniatycki13} R. Cushman and J. \'{S}niatycki,
\textquotedblleft Bohr-Sommerfeld-Heisenberg theory in geometric
quantization\textquotedblright \textit{\ J. Fixed Point Theory Appl.} 
\textbf{13} (2013) 3--24.

\bibitem{cushman-sniatycki14} R. Cushman and J. \'{S}niatycki,
\textquotedblleft On Bohr-Sommerfeld-Heisenberg
quantization, \textit{Journal of geometry and symmetry in physics} \textbf{35} (2014) 11--19.

\bibitem{dirac 25} P.A.M. Dirac, \textquotedblleft The Fundamental Equations
of Quantum Mechanics\textquotedblright , \textit{Proceedings of the Royal
Society A: Mathematical, Physical and Engineering Sciences,} \textbf{109} (752)
(1925) 642--653.%

\bibitem{dirac 30} P.A.M. Dirac, \textit{The Principles of Quantum Mechanics}, Clarendon Press, Oxford, 1930.

\bibitem{duistermaat} J.J. Duistermaat, \textquotedblleft Oscillatory
integrals, Lagrange immersions and unfolding of
singularities\textquotedblright , \textit{Comm. Pure Appl. Math.}, \textbf{27} (1974) 207--281. 

\bibitem{heisenberg} W. Heisenberg, \textquotedblleft \"{U}ber die
quantentheoretische Undeutung kinematischer und mechanischer
Beziehungen\textquotedblright , \textit{Z. f\"{u}r Phys.}, \textbf{33} (1925)
879--893.

\bibitem{hirzebruch} F. Hirzebruch, \textit{Topological methods in algebraic geometry}, 
Grundlehren der Math. Wiss. vol. \textbf{131} Springer, New York, 1966.

\bibitem{mehra-rechenberg} J. Mehra and H. Rechenberg, \textit{The
Historical Development of Quantum Theory}, vol. 2, Springer-Verlag, New York
1982.

\bibitem{richter-dullin-waalkens-wiersig} P. Richter, H. Dullin, H. Waalkens, and 
J. Wiersig, \textquotedblleft Spherical pendulum, actions, and spin\textquotedblright , 
\textit{J. Phys. Chem.} \textbf{100}(1996) 19124--19135.

\bibitem{schrodinger} E. Schr\"{o}dinger, \textquotedblleft Quantisierung
als Eigenwert Problem\textquotedblright , \textit{Annalen der Physik} \textbf{79} (1926) 361--378 and 489--527.

\bibitem{schrodinger2} E. Schr\"{o}dinger, \textquotedblleft On the relation
of the Heisenberg-Born-Jordan Quantum Mechanics and Mine\textquotedblright , 
\textit{Annalen der Physik}, \textbf{79} (1926) 734--756.

\bibitem{sniatycki75} J. \'{S}niatycki, \textquotedblleft Bohr-Sommerfeld
conditions in geometric quantization\textquotedblright , \textit{Rep. Math. Phys.} \textbf{7} (1975) 303--311.

\bibitem{sniatycki80} J. \'{S}niatycki, \textit{Geometric Quantization and
Quantum Mechanics}, Springer-Verlag, New York 1980.

\bibitem{von neumann} J. von Neumann, \textquotedblleft
Wahrscheinlichkeitstheoretischer Aufbau der
Quantenmechanik\textquotedblright\ \ G\"{o}ttinger Nachrichten, sitzung 11
November, (1927) 245--272.

\bibitem{vungoc} S. Vu Ngoc, \textquotedblleft Quantum monodromy in integrable systems\textquotedblright , 
\textit{Commun. Math. Phys.} \textbf{203} (1999) 465--479. 

\bibitem{ward-volkmer} D.W. Ward and S. Volkmer, \textquotedblleft How to
Derive the Schr\"{o}dinger Equation\textquotedblright ,
\texttt{arXiv:physics/0610121v1}. 

\bibitem{wentzel} G. Wentzel, \textquotedblleft Die mehrfach periodischen Systeme in der 
Quantenmechanik\textquotedblright , \textit{Z. f\"{u}r Phys.} \textbf{37} (1926) 80--94. 

\bibitem{woodhouse} N.M.J. Woodhouse, \textit{Geometric quantization, $2^{\mathrm{nd}}$ edition}, Oxford University Press, Oxford, UK, 1997.


\end{thebibliography}
\end{document}